\documentclass[preprint, 12pt]{elsarticle}
\usepackage{graphicx}
\usepackage{amsmath}
\usepackage{amsfonts}
\usepackage{amssymb}
\usepackage{amsthm}
\usepackage{subfig}
\usepackage{bm}
\usepackage{url}
\usepackage{verbatim}

\newcommand{\eps}{\varepsilon}
\renewcommand{\d}{\mathrm{d}}

\renewcommand{\vec}[1]{{\mathbf{#1}}}

\begin{document}


\title{Streamline integration as a method for two-dimensional elliptic grid generation}
\author{M.~Wiesenberger}
\ead{Matthias.Wiesenberger@uibk.ac.at}
\author{M.~Held}
\address{Institute for Ion Physics and Applied Physics,  Universit\"at 
   Innsbruck,  A-6020 Innsbruck, Austria}
\author{L.~Einkemmer}
\address{Numerical Analysis group,  Universit\"at Innsbruck, A-6020 Innsbruck, Austria}
 
\begin{abstract}
We propose a new numerical algorithm to construct a structured
numerical elliptic grid of a doubly connected domain. Our method is 
applicable to domains with boundaries defined by two contour lines of a two-dimensional function. 
The resulting grids are orthogonal to the boundary. Grid points as well as
the elements of the Jacobian matrix can be computed efficiently and up to machine precision. 
In the simplest case we construct conformal grids, yet with the help of weight functions and monitor metrics we can control the
distribution of cells across the domain. 
Our algorithm is parallelizable and easy to implement with
elementary numerical methods. 
We assess the quality of grids by considering both the distribution of cell sizes and the accuracy of the solution to 
elliptic problems. 
Among the tested grids these key properties are best fulfilled by the grid constructed with the monitor metric approach.

%
\end{abstract}

\maketitle
\section{Introduction} \label{sec:introduction}
A structured numerical grid is generated
by the (numerical) coordinate transformation of a rectangular ``computational
domain'' to the ``physical domain'' of interest~\cite{Thompson, Grid, Liseikin}. 
Compared to unstructured grids a structured grid allows for an easier and
computationally cheaper implementation of derivatives since the computational 
domain is a tensor product grid (a rectangle).
Generally, the same numerical techniques as in a Cartesian 
coordinate system can be used. 
Still, a structured grid must fulfill certain qualities in order to be useful 
for practical numerical computations. 
Among others are smooth and nonoverlapping coordinate lines, boundary orthogonality 
and a homogeneous distributions of cells across the domain. 
The latter condition is important in advection schemes. Too small 
cells deteriorate the CFL condition in explicit schemes and the 
condition number of the discretization matrix in implicit schemes, while large
cells deteriorate the accuracy. 
Boundary orthogonality is important in order to implement Neumann boundary
conditions.
As Reference~\cite{Thompson} pointed out it is 
desirable for the grid to be at least near-orthogonal to the boundary. Although 
it is possible to represent Neumann boundary conditions in a curvilinear grid, the
accuracy of the discretization deteriorates and the implementation is more involved
than in a boundary-orthogonal grid.  

In the literature on grid generation it has been established that elliptic grids 
are the ones that best fulfill these qualities~\cite{Thompson, Grid, Liseikin}.
These grids are generated by a coordinate transformation that fulfills 
an elliptic equation. 
However, the generation of these grids based on the algorithm proposed by Thompson, Thames and Mastin (TTM) is numerically involved~\cite{Thompson1977, Liseikin}. 
A special class of elliptic grids is generated by conformal mappings. 
The coordinates obey the Cauchy--Riemann differential equations
\begin{align}
\frac{\partial u}{\partial x} = \frac{\partial v }{\partial y} \quad
\frac{\partial u}{\partial y} = -\frac{\partial v }{\partial x}
\label{eq:CR}
\end{align}
They are orthogonal, smooth, and preserve the aspect ratio of the 
cells from the computational into the physical space. 
There are many methods to numerically construct a conformal mapping~\cite{Papamichael}, e.g. boundary integral element methods, however their implementation can be tricky. 
Among further improvements to the original TTM method are grid adaption methods and 
the monitor metric approach~\cite{Glasser2006, Vaseva2009, Liseikin}.
Both of these approaches modify the elliptic equation used 
to construct the coordinate system to a more suitable form. 
By choosing specific forms the distribution of grid cells is controlled, which
makes the coordinate transformation more flexible than conformal mappings.
The main numerical difficulty in these techniques is the solution of
the nonlinear inverted Beltrami equation.
This is required because a (not yet constructed) numerical grid is needed in order to 
discretize and solve the elliptic equation directly. If the boundary lines
are e.g.~given by a parametric representation, this conundrum 
cannot be easily solved. 
In this contribution we show how this problem can be avoided if the 
boundary lines of the physical domain are given by a two-dimensional 
function. 

One important application, where this is indeed the case, are tokamaks. 
These are magnetic fusion devices that use nearly 
axisymmetric magnetic fields to confine a plasma~\cite{Wesson}. 
In a two-dimensional plane of
fixed toroidal angle the magnetic field-lines, in the idealized case, lie within so-called flux-surfaces. 
These surfaces are in fact the isosurfaces of the poloidal magnetic flux $\psi(x,y)$ described by the Grad-Shafranov equation. 
A typical numerical modelling scenario involves the study of plasma 
dynamics between two flux-surfaces $\psi_0$ and $\psi_1$. 
Structured and unstructured numerical grids have been proposed 
for the numerical discretization of this region~\cite{hamada62, boozer80, boozer81, Grimm1983, cheng92, park08, Czarny2008, Ribeiro2010}. 
Let us emphasize, however, that the algorithms developed are in no way
restricted to magnetic fusion applications. In particular, the monitor
metric approach gives sufficient flexibility to handle large classes
of problems. For example, in certain pollution models a system of
diffusion-advection-reaction equations have to be solved for a given
velocity field (the velocity field is determined by analytic modeling,
measurements, or computer simulation). In this case the stream
function corresponding to the velocity field takes on the role of
$\psi$~\cite{Oliver2012}.
A popular choice for tokamak magnetic fields are so-called flux or magnetic coordinates~\cite{haeseleer}, in which the magnetic flux $\psi$ or a function of it is the first coordinate.
One particular class of such coordinate systems is known as PEST coordinates \cite{Grimm1983}. In this case, the geometric toroidal angle \(\varphi\) is the third coordinate, while a poloidal angle-like coordinate is constructed implicitly by choosing the volume form 
of the coordinate transformation. 
Other choices of flux coordinates exist, where the geometric toroidal angle \(\varphi\) differs from the toroidal flux angle \(\varphi_f\). 
For example, Boozer~\cite{boozer80,boozer81} and Hamada~\cite{hamada62} coordinates. In general, a flux aligned coordinate system is not orthogonal however.
Flux aligned grids are commonly constructed by 
integrating the streamlines of the vector field tangential to the isosurfaces of $\psi$~\cite{haeseleer, scott_2001}. 
Reference~\cite{Ribeiro2010} constructed a near-conformal coordinate system that
is aligned to the magnetic flux-surfaces in this way. 
The coordinates are near-conformal in the sense that the grid-deformation is small.
Unfortunately, the coordinate lines are not orthogonal.

The previous examples show that it is possible to construct a 
structured grid of a domain defined by the contour lines of a two-dimensional function $\psi(x,y)$ (the flux).
Now recall that this is exactly the requirement for the discretization 
of an elliptic equation on this domain. The inverted Beltrami equation 
in the TTM and related methods is no longer needed.

We therefore propose to construct an elliptic grid in three steps. 
First, a flux aligned grid is constructed by one of the methods that have been 
proposed in the literature.  
Then in a second step a suitably chosen elliptic equation is transformed to, 
discretized and solved in this coordinate system. 
We investigate the simple Laplace equation, the 
adapted equation and the 
elliptic equation with monitor metric suggested by References~\cite{Liseikin, Glasser2006, Vaseva2009}. 
Finally, we treat the solution of the elliptic
equation as a new flux function. We can therefore use an adapted version 
of the algorithm in the first step to construct a second coordinate transformation
from the $\psi$-aligned coordinates to coordinates aligned to the solution of 
the second step.
The final transformation then consists of the two consecutive coordinate
transformations from the first and the third step.

In Section~\ref{sec:geometry} we derive our algorithm by using methods
from differential geometry. We first explain the approach of streamline integration
by exemplarily constructing an orthogonal flux aligned grid. 
Then we show how to 
construct the conformal grid, the adapted grid and the grid using 
a monitor metric with our new method.
In Section~\ref{sec:numerics}
we numerically show that with our algorithm the coordinate map and its derivatives can be computed efficiently and up to machine precision.
Furthermore, we plot example grids for a domain relevant in 
magnetic confinement fusion, 
which are suitable for edge turbulence simulations. 
We compare the grids to
flux aligned grids and assess the quality of our grids with two different solutions
of an elliptic equation and the computation of maximal and minimal cell sizes. 
We conclude in Section~\ref{sec:conclusion}.

\section{High precision elliptic grid generation} \label{sec:geometry}

Given is a function $\psi(x,y)$ in Cartesian coordinates. 
We want to construct a grid on the region 
bounded by the two lines
$\psi(x,y) = \psi_0$ and $\psi(x,y)=\psi_1$ with $\psi_0\neq\psi_1$.
The derivative of the function $\psi$ may not vanish within this region and on the boundary and 
we further assume that the region is topologically a ring. 
Note here that this excludes the description of domains with an X-point or O-point. 
The numerical grid is described by a mapping of the discretization of the rectangular
computational domain $(u,v) \in [0,u_1]\times[0,2\pi]$ to the physical domain $(x,y)$. $u_1$ is an unknown, which the grid generation process has to provide. 

\begin{figure}[htbp]
\centering
\includegraphics[trim = 0px 0px 0px 0px, clip, scale=1.0]{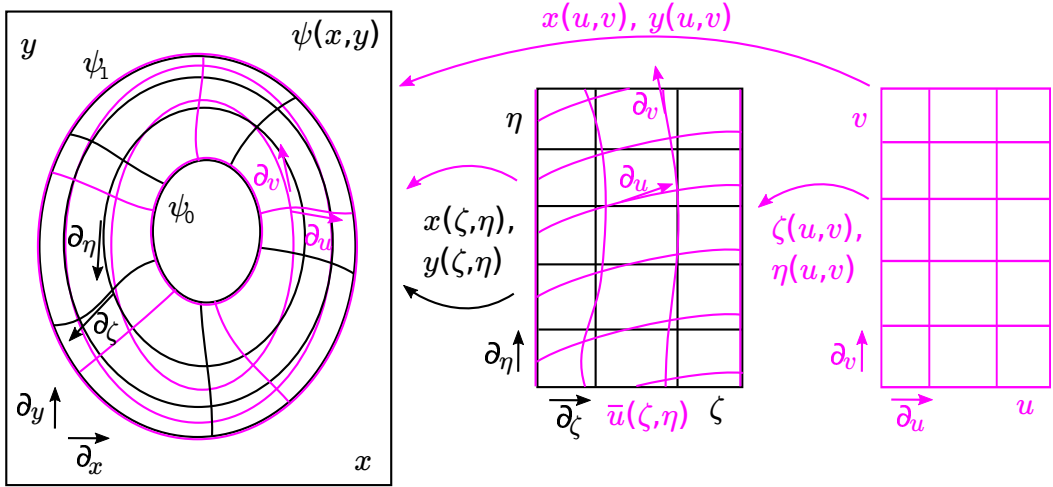}
\caption{
  Sketch of the coordinate systems, coordinate lines and basis vector fields 
  involved in our method.
}
\label{fig:sketch}
\end{figure}
As mentioned in the introduction and illustrated in Fig.~\ref{fig:sketch} the idea of our algorithm are two consecutive coordinate transformations constructed by streamline integration. 
We denote the first transformation with $x(\zeta, \eta)$, $y(\zeta, \eta)$, where
the $\zeta$ coordinate is aligned with $\psi$ and $\eta$ is an angle-like coordinate. The solution of the 
elliptic equation on this flux-aligned coordinate system is denoted with $\bar u(\zeta, \eta)$. 
The second coordinate transformation is then denoted $\zeta(u,v)$, $\eta(u,v)$ with $u$ aligned to $\bar u$.

When deriving the algorithm we use basic methods and notational 
conventions of differential geometry.  
Here, we recommend the excellent Reference~\cite{Frankel} as an introduction to the topic.  
We do so since in this approach the separate roles of the metric tensor, 
the coordinate system and its base vector fields are very clear. This is 
paramount for a concise description of our method.

In this Section we first show how the
basis vector fields can be integrated to construct coordinate lines
 with the example of orthogonal coordinates in Section~\ref{sec:orthogonal}. This follows
an introduction of some helpful quantities, our notation and streamlines
in Section~\ref{sec:preliminaries}. 
In the next step we transform, discretize and solve an elliptic equation on the flux aligned
coordinate system.
The solution $\bar u(\zeta, \eta)$ then takes the role of a new flux function in the flux aligned coordinates. 
We can therefore repeat the streamline integration method in the flux aligned
coordinate system in order to construct coordinate lines of the final $u,v$ coordinates.
We discuss three examples of this method.
In 
Section~\ref{sec:conformal} we consider the simple Laplace and the Cauchy-Riemann equations~\eqref{eq:CR}
in order to construct conformal coordinates.
In Section~\ref{sec:adaption} we modifiy the elliptic equation to allow grid adaption, before
we introduce the monitor metric in Section~\ref{sec:metric}. 
Finally, we present our algorithm in Section~\ref{sec:elliptic}. 

\subsection{Preliminaries} \label{sec:preliminaries}
In a curvilinear
coordinate system $(\zeta,\eta)$ the components of the metric tensor and its inverse
take the form 
\begin{align}
  \vec g(\zeta,\eta) = \begin{pmatrix}
    g_{\zeta\zeta} & g_{\zeta\eta}  \\
    g_{\zeta\eta} & g_{\eta\eta}  
  \end{pmatrix}
  \quad
  \vec g^{-1}(\zeta,\eta) = \begin{pmatrix}
    g^{\zeta\zeta} & g^{\zeta\eta}  \\
    g^{\zeta\eta} & g^{\eta\eta}  
  \end{pmatrix}
  \label{}
\end{align}
We denote $g = \det{\vec g} = g_{\zeta\zeta}g_{\eta\eta}-g_{\zeta\eta}^2$.
For Cartesian coordinates $(x,y)$ the elements of the inverse metric tensor are 
transformed by
\begin{subequations}
\begin{align}
  g^{\zeta\zeta} &=  \zeta_x^2+\zeta_y^2\\
  g^{\zeta\eta} &=  \zeta_x \eta_x+\zeta_y \eta_y\\
  g^{\eta\eta} &=  \eta_x^2+\eta_y^2
  \label{}
\end{align}
\end{subequations}
since $g^{xx}=g^{yy}=1$ and $g^{xy}=0$. 
Given $\zeta(x,y)$, $\eta(x,y)$ and its inverse $x(\zeta,\eta)$, $y(\zeta,\eta)$ recall that their Jacobian matrices are related by
\begin{align}
  \begin{pmatrix}
    x_\zeta & x_\eta \\
    y_\zeta & y_\eta
  \end{pmatrix} \bigg\rvert_{\zeta(x,y),\eta(x,y)}
  =
  \frac{1}{\zeta_x\eta_y - \zeta_y\eta_x }\begin{pmatrix}
    \eta_y & -\zeta_y \\
    -\eta_x & \zeta_x
  \end{pmatrix}\bigg\rvert_{x,y} 
  \label{eq:inverse}
\end{align}
With the rules of tensor transformation it is easy to prove that the element of the volume form $\sqrt{g}$ is related to the Jacobian via
\begin{align}
\sqrt{g} = (x_\zeta y_\eta-y_\zeta x_\eta) = (\eta_y\zeta_x -\eta_x\zeta_y)^{-1}
  \label{eq:vol}
\end{align}
The components of the gradient operator and the divergence in an arbitrary coordinate
system read
\begin{subequations}
\begin{align}
  \left( \nabla f \right)^i &= g^{ij}\partial_j f \\
  \nabla \cdot \vec A &= \frac{1}{\sqrt g}\partial_i \left( \sqrt{g} A^i  \right),
\end{align}
  \label{eq:arbitrary}
\end{subequations}
where we sum over repeated indices $i,j\in\{\zeta,\eta\}$ and define $\partial_\zeta \equiv \partial/\partial \zeta$, $\partial_\eta \equiv \partial/\partial \eta$.
%
Finally, we introduce the 
geometrical poloidal angle 
\begin{align}
  \theta (x,y) = \begin{cases}
    +\arccos\left( \frac{x-x_0}{\sqrt{(x-x_0)^2 + (y-y_0)^2}} \right) \text{ for } y\geq y_0 \\
    -\arccos\left( \frac{x-x_0}{\sqrt{(x-x_0)^2 + (y-y_0)^2}} \right) \text{ for } y< y_0 
  \end{cases}
  \label{eq:deftheta}
\end{align}
such that the differential 1-form
\begin{align}
  \d \theta = 
    -\frac{y-y_0}{(x-x_0)^2+(y-y_0)^2} \d x
    +\frac{x-x_0}{(x-x_0)^2+(y-y_0)^2} \d y,
  \label{eq:theta}
\end{align}
where $(x_0,y_0)$ is any point inside the region bounded by $\psi_0$.

Streamline integration is the central part of our algorithm. Recall that
given a vector field $v(x,y)=v^x(x,y)\partial_x + v^y(x,y)\partial_y$ its streamlines are given by the equation
\begin{subequations}
\begin{align}
  \frac{\d x}{\d t } = v^x(x,y)|_{x(t), y(t)} \\
  \frac{\d y}{\d t } = v^y(x,y)|_{x(t), y(t)}
\end{align}
\label{eq:streamlines}
\end{subequations}
where $t$ is a parameter. Recall here that in differential geometry 
the directional derivatives $\partial_x$ and
$\partial_y$ are the base vector fields of the coordinate system\footnote{ 
  Some textbooks (e.g.~\cite{haeseleer}) introduce the notation $\vec e_i := \partial \vec x /\partial x^i$ 
  and $\vec e^i := \nabla x^i$. 
  While this formulation is suitable in many situations we refrain from 
  using it since it mixes the metric into the basis vectors through the use of 
  the gradient (cf. Eq.~\eqref{eq:arbitrary}), which is unpractical for our purposes.}.
Now recall that if we have a function $f(x,y)$ 
such that $f(x(t), y(t))$ is a one-to-one map from $t$ to $f$ we can re-parameterize
Eq.~\eqref{eq:streamlines} by
\begin{subequations}
\begin{align}
  \frac{\d x}{\d f } = \frac{\d x/\d t|_{t(f)}}{\d f/\d t|_{t(f)}} = 
  \frac{v^x(x,y) }{ (v^x\partial_xf + v^y\partial_yf)(x,y) }\bigg |_{x(f), y(f)}\\
  \frac{\d y}{\d f } = \frac{\d y/\d t|_{t(f)}}{\d f/\d t|_{t(f)}} = 
  \frac{v^y(x,y) }{ (v^x\partial_xf + v^y\partial_yf)(x,y) }\bigg |_{x(f), y(f)}
\end{align}
\label{eq:reparameter}
\end{subequations}
Furthermore, the derivative of any function $g(x,y)$ along the streamlines of $v$ 
parameterized by $f$ reads
\begin{align}
  \frac{\d g}{\d f }\bigg |_{x(f),y(f)} = 
  \frac{(v^x\partial_x g + v^y\partial_y g)(x,y) }{ (v^x\partial_xf + v^y\partial_yf)(x,y) }\bigg |_{x(f),y(f)}
  \label{}
\end{align}

Figuratively, in any coordinate system $(\zeta,\eta)$ the 1-forms $\d \zeta$ and $\d \eta$ (the contravariant basis) are visualized by the lines (surfaces in higher dimensions) 
of constant $\zeta$ and $\eta$. 
At the same time the streamlines of the vector fields $\partial_\zeta$ and $\partial_\eta$ (the covariant basis) 
give the coordinate lines of $\zeta$ and $\eta$ (cf. Fig.~\ref{fig:sketch}). 
For example, if we hold $\eta$ constant and vary $\zeta$, we go along a streamline of 
$\partial_\zeta$. This implies that in two dimensions $\d \eta$ and $\partial_\zeta$ 
trace the same line. The vector fields $\nabla \zeta$ and $\nabla \eta$ are associated to 
$\d \zeta$ and $\d \eta$ through the metric tensor by Eq.~\eqref{eq:arbitrary} and are the vector fields that
are everywhere perpendicular to the lines of constant $\zeta$ and $\eta$, respectively. 
It is important to realize that in general curvilinear coordinates the vector
fields $\nabla \zeta$ and $\nabla \eta$ point in 
different directions than $\partial_\zeta$ and $\partial_\eta$.
The central point in our algorithm is the realization that once we can 
express $\partial_\zeta$ and $\partial_\eta$ in terms of $\partial_x$ and 
$\partial_y$ we can immediately construct the coordinate transformation 
by integrating 
streamlines of $\partial_\zeta$ and $\partial_\eta$ using Eq.~\eqref{eq:streamlines}. This 
holds true even if $(x,y)$ were curvilinear coordinates. The components of the 
one-forms $\d \zeta$ and $\d\eta$ in terms of $\d x$ and $\d y$ form 
the elements of the Jacobian matrix of the transformation. These 
are also necessary in order to transform any tensor (including the metric) 
from the old to the new coordinate system.

\subsection{Orthogonal coordinates} \label{sec:orthogonal}
In general, orthogonal coordinates $\zeta, \eta$ with $\zeta$ aligned to $\psi$ are described by
\begin{subequations}
\begin{align}
  \d \zeta & = \zeta_x\d x +\zeta_y\d y =  f(\psi)(\psi_x \d x + \psi_y \d y) \\
  \d \eta  & = \eta_x \d x + \eta_y \d y = h(x,y) ( -\psi_y \d x + \psi_x \d y)
\end{align}
  \label{eq:orthogonal}
\end{subequations}
With Eq.~\eqref{eq:orthogonal} we have $g^{\zeta\eta} = \zeta_x\eta_x + \zeta_y\eta_y = 0$,
$g^{\zeta\zeta} = (\nabla\psi)^2 f^2$, $g^{\eta\eta} = (\nabla\psi)^2h^2$ and $\sqrt{g}^{\,-1} = (\nabla\psi)^2 h f$. 
From Eq.~\eqref{eq:inverse} we directly see that the basis vector fields are
\begin{subequations}
\begin{align}
  \partial_\zeta&= x_\zeta\partial_x + y_\zeta\partial_y = \frac{1}{(\nabla\psi)^2f} (\psi_x \partial_x + \psi_y\partial_y) \\
  \partial_\eta &=  x_\eta\partial_x + y_\eta\partial_y =  \frac{1}{(\nabla\psi)^2h} (-\psi_y \partial_x + \psi_x\partial_y) 
\label{eq:orthogonal_linesb}
\end{align}
\label{eq:orthogonal_lines}
\end{subequations}
i.e. $\partial_\zeta$ points into the direction of the gradient of $\psi$ and $\partial_\eta$ into the direction of constant $\psi=\text{const}$ surfaces.
Now, the coordinate system is defined up to the functions $f(\psi)$ and $h(x,y)$.
Note that $f$ must be a function of $\psi$ only since the restriction $\d(\d \zeta)=0$ must hold.
Furthermore, $f(\psi) = \d\zeta/ \d\psi \neq 0 $ is in principle an arbitrary function, yet we choose $f(\psi) = f_0 = \text{const}$.
With this choice we directly get
\begin{align}
  \zeta(x,y) = f_0(\psi(x,y)-\psi_0)
  \label{eq:zeta}
\end{align}
Note that $\zeta_0=0$ and $\zeta_1=f_0(\psi_1-\psi_0)$.
The up to now undefined function $h(x,y)$ is not arbitrary since $\d(\d \eta) = 0$ must hold. This is the requirement that $\d \eta$ must be
a closed form in order for the potential $\eta$ to exist. 
We can express this as
\begin{align}
  (\psi_x\partial_x + \psi_y\partial_y) h = f(\nabla\psi)^2 \partial_\zeta h= -h\Delta \psi 
  \label{eq:hequation}
\end{align}
where $\Delta \psi = \psi_{xx} + \psi_{yy}$ is the two-dimensional Laplacian.
Let us remark here that Eq.~\eqref{eq:hequation} can be written as $\nabla\cdot\left( h\nabla\psi \right)=0$, which makes the orthogonal grid an elliptic grid
with adaption function $h$ as becomes evident later. 
In order to integrate this equation we need an initial condition for $h$. 
We choose to first discretize the line given by $\psi(x,y) = \psi_0$. 

As already mentioned $\partial_\eta$ is the vector field the streamlines of which give the 
coordinate lines for $\eta$. We choose $h(x,y) = \text{const}$ on $\psi_0$. To this end
we parameterize the coordinate line by $\theta$ (cf. Eq.~\eqref{eq:reparameter})
\begin{subequations}
\begin{align}
  \left . \frac{\d x}{\d \theta}\right|_{\zeta=0}=\frac{x_\eta}{\theta_\eta} = \frac{-\psi_y}{\psi_x\theta_y - \psi_y \theta_x}\\
  \left . \frac{\d y}{\d \theta}\right|_{\zeta=0}=\frac{y_\eta}{\theta_\eta} = \frac{\psi_x}{\psi_x\theta_y - \psi_y \theta_x}\\
  \left . \frac{\d \eta}{\d \theta}\right|_{\zeta=0}=\frac{1}{\theta_\eta} = \frac{(\nabla\psi)^2h(\psi_0)}{\psi_x\theta_y - \psi_y \theta_x}
\end{align}
  \label{eq:etaline}
\end{subequations}
Let us define $h(\psi_0)$ such that $\eta \in [0,2\pi]$, that is,
\[ 2\pi = \oint_{\psi=\psi_0}\d \eta = \oint_0^{2\pi} \frac{\d \eta}{\d\theta}\bigg|_{\zeta=0}\d \theta \nonumber \]
  or 
\begin{align}
  f_0 := h(\psi_0) = \frac{2\pi}{\int_0^{2\pi}\d\theta \frac{(\nabla\psi)^2}{\psi_x\theta_y - \psi_y \theta_x}}
  \label{eq:definef}
\end{align}
Here, we also fixed the constant $f_0$ such that our coordinate system $\zeta, \eta$ fulfills the Cauchy--Riemann condition on the boundary line $\psi_0$. 
As initial point for the integration of Eq.~\eqref{eq:etaline} we can use any point with $\psi(x,y) = \psi_0$.
We then use $h(\psi_0)$ on the flux-surface $\psi_0$ as initial condition for the integration 
of Eq.~\eqref{eq:hequation}.

We obtain coordinate lines by integrating the vector fields $\partial_\zeta$ and $\partial_\eta$ given in \eqref{eq:orthogonal_lines}. We start the construction by integrating $\partial_\eta$ for $\psi = \psi_0$, i.e. $\zeta=0$. This can be done since $h|_{\psi_0}$ is known. 
The obtained points serve as starting points for the integration of $\partial_\zeta = f_0^{-1} \partial_\psi$.
In order to get $h$ we simply integrate Eq.~\eqref{eq:hequation}
\begin{subequations}
\begin{align}
    \left .\frac{\d x}{\d \zeta}\right |_{\eta=\text{const}} &= \frac{\psi_x}{f_0(\nabla\psi)^2}\\
    \left .\frac{\d y}{\d \zeta}\right|_{\eta=\text{const}} &= \frac{\psi_y}{f_0(\nabla\psi)^2}\\
    \left .\frac{\d h}{\d \zeta}\right|_{\eta=\text{const}} &= - \frac{\Delta \psi}{f_0(\nabla\psi)^2} h
\end{align}
  \label{eq:etacoordinates}
\end{subequations}
Note that if $\Delta\psi=0$, we directly get a conformal grid with this algorithm. This can be seen as then $h(\zeta,\eta) = f_0$. In~\ref{app:conf-field}
we briefly study the class of functions $\psi$ that are solutions of the Grad--Shafranov equation and satisfy $\Delta \psi=0$. This, however, is not true in general.

Let us further remark on the sign of $f_0$. It is our goal to construct a right handed
coordinate system and to have $\zeta_1>\zeta_0=0$.
The curves of constant $\theta$ surround $x_0, y_0$ in a mathematically positive direction. 
That means that Eq.~\eqref{eq:definef} implies that $f_0>0$ if $\nabla \psi$ points away from $x_0, y_0$. If this is not the case, we obtain $f_0 < 0$. 
On the other hand if $\zeta$ should increase from $\psi_0$ to $\psi_1$, we
need $f_0 <0$ for $\psi_1<\psi_0$ and $f_0>0$ for $\psi_1>\psi_0$. 
We thus simply take the absolute value of Eq.~\eqref{eq:definef} and multiply 
by $-1$ if $\psi_1<\psi_0$. 
For ease of notation we do so also in the following without further notice.

Let us finally summarize the grid generation in the following algorithm; we 
assume that the $\zeta$ space is discretized by a list of 
not necessarily equidistant values $\zeta_i$ with $i = 0,1,\dots N_\zeta-1$ 
and $\eta$ is discretized by a list of $\eta_j$ with $j = 0,1,\dots N_\eta-1$:
\begin{enumerate}
  \item Find an arbitrary point $(x,y)$ with $\psi(x,y) = \psi_0$
    and a point $x_0, y_0$ within the region bound by $\psi(x,y) = \psi_0$ for the definition of $\theta$ in Eq.~\eqref{eq:deftheta}
  \item Integrate Eq.~\eqref{eq:etaline} with $h=1$ over $\Theta=[0,2\pi]$ and use Eq.~\eqref{eq:definef} to compute $f \equiv f_0$ and $h(\psi_0)$.
    Use any convenient method for the integration of ordinary differential equations.
  \item Integrate one streamline of Eq.~\eqref{eq:orthogonal_linesb} with $h=f_0$
    from $\eta = 0..\eta_j$ for all $j$.
    The result is a list of $N_\eta$ coordinates $x_j, y_j$ on the $\psi_0$ surface.
  \item Using this list as starting values integrate Eq.~\eqref{eq:etacoordinates}
    from $\zeta=0\dots\zeta_i$ for all $i$. This gives the map $x(\zeta_i, \eta_j), y(\zeta_i, \eta_j)$ for all $i$ and $j$.
  \item Last, using the resulting list of coordinates and 
    Eq.~\eqref{eq:orthogonal} evaluate the derivatives 
    $\zeta_x(\zeta_i,\eta_j)$, $\zeta_y(\zeta_i, \eta_j)$, $\eta_x(\zeta_i,\eta_j)$, and $\eta_y(\zeta_i, \eta_j)$ for all $i$ and $j$.
\end{enumerate}

\subsection{Conformal coordinates} \label{sec:conformal}
A conformal mapping $u(x,y), v(x,y)$ has to satisfy the Cauchy--Riemann equations given in Eq.~\eqref{eq:CR}.
A direct consequence is that $u$ and $v$ are harmonic functions 
\begin{align}
  \Delta u = \Delta v = 0
  \label{eq:harmonic}
\end{align}
Here, $\Delta=\nabla^2$ is the two-dimensional Laplacian with the divergence
and gradient operators defined in~\eqref{eq:arbitrary}.
First, we note that Eq.~\eqref{eq:harmonic} holds in every coordinate system.
Let us assume that we have constructed flux aligned  coordinates $(\zeta, \eta)$. These can be, but not necessarily have to be,
the orthogonal coordinates introduced in the last section. 
Now, in order to construct conformal coordinates $u$, $v$ we first define 
\begin{align}
u(\zeta,\eta):=c_0(\bar u(\zeta,\eta)-\psi_0)
  \label{eq:ubar}
\end{align}
and thus
\begin{align}
  \Delta \bar u(\zeta,\eta) = 0 
  \label{eq:ubar_harmonic}
\end{align}
where $\bar u( 0, \eta) = \psi_0$ and $\bar u( \zeta_1, \eta) =\psi_1$ fulfills Dirichlet boundary conditions in $\zeta$. 
In $\eta$ we have periodic boundary conditions. 

Note the analogy between Eq.~\eqref{eq:ubar} and Eq.~\eqref{eq:zeta}. Now, $\bar u$ is
equal to $\psi$ at the boundaries and its Laplacian vanishes in between. 
In fact, $\bar u$ takes the role of $\psi$ in the following coordinate
transformation. 
We introduce $c_0$ as a normalization constant with the same role as $f_0$ in the orthogonal coordinate transformation. 
Having $\bar u(\zeta, \eta)$, 
our idea is to construct the basis one-forms $\d u$ and $\d v$  in terms 
of $\d \zeta$ and $\d \eta$ by transforming the Cauchy-Riemann equations to the
$\zeta, \eta$ coordinate system. Analogues to the algorithm in Section~\ref{sec:orthogonal} we can then construct the basis vector fields $\partial_u$ and $\partial_v$, appropriately choose a normalization and then use streamline
integration in the $\zeta, \eta$ coordinate system to construct the coordinates. 
From the basis one-forms $\d u$ and $\d v$ we get the elements of the Jacobian 
matrix of the transformation. 
Before we do this in detail however, let us first discuss some alternative elliptic
equations to the simple Eq.~\eqref{eq:harmonic}.

\subsection{Grid adaption} \label{sec:adaption}
Although the conformal grid is advantageous for elliptic equations (due to the vanishing metric
coefficients) the cell distribution is not very flexible; once the boundary is set
the conformal map is unique. We therefore have little control over the 
distribution of cells. 
We can use grid adaption techniques to overcome this restriction. 
The idea is to modify the elliptic equations that $u$ and $v$ have 
to fulfill. That is, we choose 
\begin{align}
\nabla\cdot\left( \frac{\nabla u}{w}\right) = \nabla\cdot\left( w \nabla v\right) = 0
\label{eq:adaption}
\end{align}
where $w$ is an appropriately chosen weight function. The 
cell size will be small in regions where $w$ is large and spread out in regions
where $w$ is small.
The Cauchy--Riemann equations~\eqref{eq:CR} are changed accordingly to 
\begin{align}
v_x = -\frac{u_y}{w}\quad v_y = \frac{u_x}{w}
\label{eq:CR_adaption}
\end{align}
Let us remark here that it is straightforward to implement the weight function in the orthogonal grid generation. 
In Section~\ref{sec:orthogonal} we simply replace the function $h$ by $h/w$.
Then we have
\begin{subequations}
\begin{align}
\partial_\zeta &= \frac{1}{f_0(\nabla\psi)^2} (\psi_x\partial_x + \psi_y \partial_y)\\
\partial_\eta &= \frac{w}{h(\nabla\psi)^2} (-\psi_y\partial_x + \psi_x \partial_y)
\end{align}
\end{subequations}
and $\nabla \psi\cdot \nabla (h/w) = -h/w \Delta\psi$.
A suitable choice for $w$ is 
\begin{align}
w = |\nabla\psi|
\label{eq:weight_adaption}
\end{align}
as then the angle-like coordinate $\eta$ becomes the arc length on the $\psi_0$ line. 


\subsection{Monitor metric and the heat conduction tensor} \label{sec:metric}
We follow Reference~\cite{Liseikin, Glasser2006, Vaseva2009} and replace the canonical metric tensor $\vec g$ by
a specifically tailored tensor $\vec G$ that takes the form
\begin{align}
\vec G(x,y) = \vec T \vec T + k^2\vec N \vec N + \eps(x,y)\vec I
\label{eq:metric}
\end{align}
with $\vec T = (-\psi_y, \psi_x)$ and $\vec N = -(\psi_x, \psi_y)$. 
The vector $\vec T$ is tangential to the contour lines of $\psi$ while $\vec N$ is 
normal to it. 
We have 
\begin{align}
  \sqrt{G} = \left[(\eps+k^2(\psi_x^2+\psi_y^2))(\eps+(\psi_x^2+\psi_y^2))\right]^{-1/2}
  \label{}
\end{align}
where $k<1$ is a constant and $\eps$ is a function that is nonzero in the neighborhood
of singularities, i.e.~where $\nabla \psi(x,y)=0$.
In our work we choose $k=0.1$ and $\eps(x,y)\equiv\eps =0.001$.
Note that the scalar product induced by $\vec G$ conserves perpendicularity 
with respect to $\vec T$, i.e.~if
and only if
 any vector $\vec v \perp \vec T$ in the canonical metric, then it is also perpendicular to $\vec T$ in $\vec G$.
In our application this is important at the boundary.\\
Now, we consider the elliptic equation
\begin{align}
  \nabla\cdot(\sqrt{G}\vec G \nabla u) = 0\nonumber\\
  \partial_x(\sqrt{G} (G^{xx}\partial_x u + G^{xy}\partial_y u )) + \partial_y(\sqrt{G}(G^{yx}\partial_x u + G^{yy}\partial_y u )) = 0
  \label{eq:elliptic}
\end{align}
with Dirichlet boundary conditions. 
The resulting grid coordinate $u$ is almost perfectly aligned in regions
far away from singularities and breaks the alignment in regions where $|\nabla\psi|$ is small or vanishes.

We now take a slightly more general approach and rewrite Eq.~\eqref{eq:elliptic}
\begin{align}
  \nabla\cdot(\vec \chi \nabla u) = 0
  \label{eq:elliptic_mod}
\end{align}
where $\chi(x,y)$ is a symmetric positive-definite contravariant tensor. 
Then the conformal grid from Section~\ref{sec:conformal}, the grid adaption from Section~\ref{sec:adaption} as 
well as the monitor metric can be considered special cases. 
The grid adaption is recovered by setting $\vec \chi = 1/w\vec I$,
while the monitor metric is simply $\vec \chi = \sqrt{G}\vec G$. 
The true conformal case is, of course, recovered by setting $\vec \chi = \vec I$. 

This allows us to provide a commonly known physical interpretation of Eq.~\eqref{eq:elliptic_mod}.
If $u$ is a temperature, then the Dirichlet boundary condition 
fixes a temperature at the boundary of our domain. The tensor $\vec \chi(x,y)$ is 
then the anisotropic heat conduction tensor and the coordinate lines for $u$ are 
the isothermal lines of the steady state solution to the heat diffusion problem. 
This interpretation allows us to intuitively estimate how the
coordinate lines look like when a specific $\vec\chi$ is chosen. 
In the case of grid adaption, if the weight function is large, the heat conduction 
is low resulting in small temperature gradients and thus closely spaced grid cells. 
On the other hand, let us reconsider Eq.~\eqref{eq:metric} for the case $\eps=0$. 
Then we can write 
\begin{align}
\vec \chi = \frac{1}{k}\hat t\hat t + k \hat n\hat n 
= \chi_\parallel \hat t \hat t + \chi_\perp \hat n \hat n
  \label{eq:heat_conduction}
\end{align}
which for $k<1$ simply means that the heat conduction parallel to the magnetic field
is far stronger than perpendicular to it, which is in fact the case in an actual 
fusion reactor. 
The coordinate lines will thus tend to align with the magnetic flux
surfaces with the degree of alignment given by $k$ resulting 
in an almost aligned grid. 
In the limit of vanishing $k$ the alignment should be perfect. 
If $|\nabla\psi|$ vanishes, the tensor~\eqref{eq:metric} reduces to $\vec \chi = \vec I$.

\subsection{The elliptic grids} \label{sec:elliptic}
Suppose that we have constructed a boundary aligned grid $(\zeta,\eta)$, which
may but not necessarily has to be the orthogonal grid from Section~\ref{sec:orthogonal}. 
Now, we solve the general elliptic equation
\begin{align}
 \nabla\cdot(\vec \chi \nabla \bar u) = \partial_i(\sqrt{g} \chi^{ij}\partial_j \bar u) = 0
  \label{eq:general_elliptic}
\end{align}
in the transformed coordinate system with boundary conditions $\bar u|_{\partial\Omega}=\psi$. We set $u=c_0(\bar u -\psi_0)$. This means that we have to 
transform the conduction tensor $\chi$ from Cartesian to flux coordinates,
which is done by the well known rules of tensor transformation
\begin{subequations}
\begin{align}
  \chi^{\zeta\zeta}(\zeta, \eta) &= (\zeta_x\zeta_x \chi^{xx} + 2\zeta_x\zeta_y\chi^{xy} + \zeta_y\zeta_y \chi^{yy})|_{x(\zeta,\eta), y(\zeta,\eta)}\\
  \chi^{\zeta\eta}(\zeta, \eta) &= (\zeta_x\eta_x \chi^{xx} + (\zeta_x\eta_y + \eta_x\zeta_y)\chi^{xy} + \zeta_y\eta_y \chi^{yy})|_{x(\zeta,\eta), y(\zeta,\eta)}\\
  \chi^{\eta\eta} (\zeta, \eta) &= (\eta_x\eta_x \chi^{xx} + 2\eta_x\eta_y\chi^{xy} + \eta_y\eta_y \chi^{yy}))|_{x(\zeta,\eta), y(\zeta,\eta)}
\end{align}
\label{eq:transformationChi}
\end{subequations}
The equivalent of the Cauchy--Riemann equations in this formulation reads
\begin{subequations}
\begin{align}
  v_\zeta = -\sqrt{g}(\chi^{\eta\zeta}u_\zeta + \chi^{\eta\eta}u_\eta)\\
  v_\eta = +\sqrt{g}(\chi^{\zeta\zeta}u_\zeta + \chi^{\zeta\eta}u_\eta)
\end{align}
  \label{eq:hodge_dual}
\end{subequations}
These are constructed such that
 $\nabla\cdot((\vec \chi/ \det \chi) \nabla v) = 0$.
The interested reader might notice that Eq.~\eqref{eq:hodge_dual} just 
defines the components of the Hodge dual $\d v = \star\d u$  
if $\chi$  is interpreted as a metric.
We note that these equations are now valid for any grid that  
we use to solve Eq.~\eqref{eq:general_elliptic}. If we find a boundary
aligned grid analytically, we can start the grid construction 
directly with the solution of Eq.~\eqref{eq:general_elliptic} 
and then proceed with the conformal grid generation. 

The relevant equations for the streamline integration now read
\begin{subequations}
\begin{align}
    \d u &= c_0(\bar u_\zeta\d \zeta + \bar u_\eta \d \eta) \\
    \d v &= c_0\sqrt{g}(-(\chi^{\eta\zeta}\bar u_\zeta + \chi^{\eta\eta}\bar u_\eta)\d \zeta + (\chi^{\zeta\zeta}\bar u_\zeta + \chi^{\zeta\eta}\bar u_\eta)\d \eta)
\end{align}
\label{eq:contravariant_monitor}
\end{subequations} 
which just means that $v$
is orthogonal to $u$ in the scalar product 
generated by the symmetric tensor $\vec \chi$, which we denote by $\langle . , .\rangle$.
We have 
$J=c_0^2 \sqrt{g} \langle \nabla \bar u,\nabla \bar u \rangle
  := c_0^2\sqrt{g} (\bar u_\zeta^2\chi^{\zeta\zeta} + 2 \bar u_\zeta \bar u_\eta \chi^{\zeta\eta} + \bar u_\eta^2\chi^{\eta\eta})$ 
and 
\begin{subequations}
\begin{align}
  \partial_u &= \frac{1}{c_0\langle\nabla \bar u,\nabla \bar u\rangle} 
        (\chi^{\zeta\zeta}\bar u_\zeta + \chi^{\zeta\eta}\bar u_\eta)\partial_\zeta + 
        (\chi^{\eta\zeta}\bar u_\zeta + \chi^{\eta\eta}\bar u_\eta)\partial_\eta\label{eq:basis_monitora}\\
  \partial_v &= \frac{1}{c_0\sqrt{g}\langle\nabla \bar u,\nabla \bar u\rangle} 
          (-\bar u_\eta \partial_\zeta + \bar u_\zeta\partial_\eta)
  \label{eq:basis_monitorb}
\end{align}
\label{eq:basis_monitor}
\end{subequations} 
As for the orthogonal coordinates we have to integrate these two vector fields to construct our coordinates. 
We begin with the integration of $\partial_v$ along the $\zeta=0$ line. 
It is important to note that $\bar u_\eta|_{\zeta=0} = 0$ and thus
\begin{align}
  \eta_v(0, \eta) = \left( c_0\sqrt{g}\bar u_\zeta \chi^{\zeta\zeta} \right)|_{\zeta= 0, \eta}
  \label{eq:normalize}
\end{align}
We can use Eq.~\eqref{eq:normalize} to define $c_0$ such that $v\in[0,2\pi]$. 
In order to do so we simply integrate 
\begin{align}
  v_1 = \int_0^{2\pi} \left .\frac{\d v}{\d \eta}\right|_{\zeta=0} \d \eta = \int_0^{2\pi} \frac{1}{\eta_v}\d \eta = c_0 \int_0^{2\pi} \sqrt{g}\chi^{\zeta\zeta}\bar u_\zeta(0,\eta) \d \eta := 2\pi
  \label{eq:computec0}
\end{align}
 with $v_0 = 0$ and choose $c_0$ such that $v_1=2\pi$. This is the analogues
 equation to Eq.~\eqref{eq:definef}, with the difference that we can integrate
 Eq.~\eqref{eq:computec0} directly using numerical quadrature. 
Having done this we integrate, analog to Section~\ref{sec:orthogonal}, 
$\partial_v $ on $\zeta=0$ to get starting points for the integration of $\partial_u$ 
from $\bar u_0 = 0$ to $\bar u_1= c_0 \zeta_1$.
Note, that we can compute the components of $\d u$ and $\d v$ in terms of $\d x$ and $\d y$ by using the transformation
\begin{subequations}
\begin{align}
  u_x(\zeta, \eta)  &= u_\zeta \zeta_x + u_\eta\eta_x,\ u_y(\zeta, \eta) = u_\zeta \zeta_y + u_\eta\eta_y\\
  v_x(\zeta, \eta)  &= v_\zeta \zeta_x + v_\eta\eta_x,\ v_y(\zeta, \eta) = v_\zeta \zeta_y + v_\eta\eta_y
\end{align} 
  \label{eq:chain_rule}
\end{subequations}
as soon as the constant $c_0$ becomes available. 

Note that the resulting grid will, in general, not be orthogonal 
in the Euclidean metric. However, for all the cases we discuss
in this paper, the grid is orthogonal at the boundary.

The final algorithm now reads, assuming that $u$ is discretized
by not necessarily equidistant points $u_i$ for $i=0,1,\dots N_u-1$ and 
$v$ by $v_j$ with $j=0,1,\dots N_v-1$ and that $x(\zeta,\eta)$, $y(\zeta, \eta)$
as well as the components of the Jacobian $\zeta_x(\zeta,\eta)$, $\zeta_y(\zeta,\eta)$, $\eta_x(\zeta,\eta)$ and $\eta_y(\zeta,\eta)$ are available from the first coordinate transformation (e.g. Section~\ref{sec:orthogonal}):
\begin{enumerate}
  \item Choose either $\vec \chi = \vec I$, $\vec \chi=1/w\vec I$ or $\vec \chi = \sqrt{G}\vec G$ depening on whether a conformal, adapted or monitor grid is desired.
  \item Discretize and solve the elliptic equation~\eqref{eq:general_elliptic} on the $\zeta,\eta$ grid for $\bar u(\zeta, \eta)$ with any
    method that converges.  
    Use Eq.~\eqref{eq:transformationChi} to transform $\chi$ from Cartesian 
    to the $\zeta, \eta$ coordinate system.
    The chosen resolution
    determines the accuracy of the subsequent streamline integration.
  \item Numerically compute the derivatives $\bar u_\zeta$ and $\bar u_\eta$ 
    and construct $\eta_v^{-1}(0, \eta)$ using Eq.~\eqref{eq:normalize} for $c_0=1$
  \item Integrate Eq.~\eqref{eq:computec0} to determine $c_0$.
  \item On the $\zeta, \eta$ grid compute $\eta_v$, $\zeta_u$ and $\zeta_v$ according to Eq.~\eqref{eq:basis_monitor} as well as $u_\zeta$, $u_\eta$, $v_\zeta$ and
    $v_\eta$ according to Eq.~\eqref{eq:contravariant_monitor}. Use Eq.~\eqref{eq:chain_rule} and the Jacobian of the $\zeta, \eta$ coordinates to compute $u_x$, $u_y$, $v_x$ and $v_y$.
  \item Integrate the streamline of $\partial_v$ for $\zeta=0$ 
    from $v=0\dots v_j$ for all $j$ on the $\zeta, \eta$ grid using the normalized component $\eta_v$. 
    Interpolate $\eta_v(\zeta, \eta)$ when necessary.
  \item Using the resulting points as start values integrate $\partial_u$ from 
    $u=0\dots u_i$ for all $i$ and all points. The result is the list
    of coordinates $\zeta(u_i, v_j)$, $\eta(u_i, v_j)$ for all $i$ and $j$.
  \item Interpolate $x\left( \zeta,\eta \right)$, $y\left( \zeta, \eta \right)$, 
    $u_x(\zeta, \eta)$, $u_y(\zeta, \eta)$, $v_x(\zeta, \eta)$ and $v_y(\zeta, \eta)$ on this list. 
\end{enumerate}
There are two differences in this algorithm from the one presented in 
Section~\ref{sec:orthogonal}:
For the integration of the vector fields Eq.~\eqref{eq:basis_monitor} and the evaluation of derivatives Eq.~\eqref{eq:contravariant_monitor} we need to evaluate its components at arbitrary points.
An interpolation method is thus needed. 
In order to avoid out-of-bound errors we can artificially make the $\zeta,\eta$ box periodic.
Second, the existence of the coordinate $v$ is guaranteed by the Cauchy--Riemann equations
and thus the $h$ function does not appear. 

A suitable test for the implementation is e.g. the volume/area of the domain. 
Being an invariant
the volume must be the same regardless of the coordinate system in use. 


\section{Numerical tests \label{sec:numerics}}
We extend our numerical library FELTOR ({\url{www.github.com/feltor-dev/feltor}}) \cite{FELTORv3.1} 
with a geometry package that can handle the various grids discussed in Section~\ref{sec:geometry}. 
For our implementation we choose high order explicit Runge--Kutta methods (see, for example, \cite{hairer93}) to integrate the necessary ordinary differential equations
and discontinuous Galerkin (dG) methods~\cite{Cockburn2001, Held2016} for
the elliptic equations and the interpolation.
The grid points and metric elements are thus available up to machine precision. 
Note here that each point of the grids we discussed
in Section~\ref{sec:geometry} can be computed independently from each 
other, which provides a trivial parallelization option. 

We perform the computationally intensive inversion of the elliptic 
equation~\eqref{eq:general_elliptic} on a GPU with a conjugate gradient method.
On a single CPU a direct solver might be 
preferable but the assembly of the elliptic operator is not straightforward due
to the presence of the metric elements. 
The advantage in an iterative method is that
we do not need to assemble the whole matrix 
we only need to implement the
application to a vector.

As already discussed in the introduction in all structured grids derivatives are pulled back to a rectangular grid (the computational space).
The computational space is a product space, i.e.~the discretization 
of the derivatives is one-dimensional (cf. \cite{Einkemmer2014, Held2016}). 
This reduces the stencil of the 
discretization compared to a discretization using an unstructured two-dimensional grid and simplifies the communication pattern in a parallel implementation. 
Furthermore, the matrices have a block-diagonal form with two side bands, 
where each block has the size $P\times P$ with $P$ the number of polynomial 
coefficients in each cell.
Apart from the corner entries all blocks on the diagonals are equal. 
This further reduces the storage requirements and
in our experience increases the performance of matrix-vector multiplications by 
a factor of $2-3$ compared to a sparse matrix format (like the compressed sparse row format) that stores all non-zero elements explicitly.


In the following we will use the function $\psi$ that is given 
by Reference~\cite{Cerfon2010} as an analytical solution to the Grad-Shafranov equation. This
function is used in practice for turbulence simulations in tokamaks. 
The values of the $12$ coefficients in $\psi$ can be found in \ref{app:coeff-psi}. Furthermore, 
all programs used for this work can be found in the latest
release of FELTOR~\cite{FELTORv3.1}.
\subsection{Convergence of $\bar u$}
We first check that the numerical solution of Eq.~\eqref{eq:general_elliptic} converges as expected.  
As we do not have an analytical solution we compute the difference in the $L2$ norm of the computational space 
between two consecutive solutions $\bar u_{\text{num}}^1$ and $\bar u_{\text{num}}^2$. 
We thus define the numerical error 
\begin{align}
  \eps_{\bar u} = \left(\frac{\int_{\zeta_0}^{\zeta_1} d \zeta\int_0^{2\pi}\d \eta (\bar u_{\text{num}}^1 - \bar u_{\text{num}}^2)^2}{ \int_{\zeta_0}^{\zeta_1} d \zeta\int_0^{2\pi}\d \eta(\bar u_{\text{num}}^1)^2}\right)^{1/2}.
\label{eq:error_1}
\end{align}
\begin{table}[htbp]
\begin{center}
\begin{tabular}{|c|c|c|c|c|c|c|}
\hline
$N_\eta=10 N_\zeta$& $P=9$    & $P=11$     & $P=13$    & $P=15$   \\ \hline
         20     &  -       &   -        &       -   &  -       \\ \hline
         40 & 5.98E-05 & 8.69E-06 & 5.86E-06 & 2.87E-07 \\ \hline
         80 & 3.48E-06 & 3.40E-07 & 3.75E-08 & 1.89E-08 \\ \hline
         160 & 4.80E-08 & 3.58E-09 & 2.89E-10 & 2.64E-11 \\ \hline
         320 & 7.43E-10 & 3.02E-11 & 1.71E-12 & 3.07E-12  \\ \hline
\end{tabular}
\end{center}
\caption{ Comparison of the $L2$ errors \eqref{eq:error_1} for the solution of the elliptic equation with monitor metric~\eqref{eq:general_elliptic} on flux grid. 
    The number of cells is denoted by $N_{\eta}$ and $N_{\zeta}$, respectively, and we employ polynomials of degree $P-1$ in each cell.
}
\label{tab:convergence_u}
\end{table}

In Table~\ref{tab:convergence_u} we show the results for different high order polynomials to test the
convergence of the general elliptic equation with monitor metric~\eqref{eq:general_elliptic}. 
The region is bounded by $\psi_0 = -20$ and $\psi_1=-1$.
We use a fluxgrid with the angle-like coordinate defined as the arc length to discretize this region~\cite{Grimm1983}. The advantage 
is that this grid is quickly generated and has,
other than the orthogonal grids discussed in Section~\ref{sec:orthogonal}, 
a fairly homogeneous distribution 
of cells throughout the domain. As we use Dirichlet boundary conditions 
the nonorthogonality is not an issue. 
We observe a quick convergence until a relative error of around $10^{-12}$ is reached. This apparent upper bound is due to the accuracy of the residuum in the 
conjugate gradient solver at $10^{-11}$.
Let us note that we observe similar convergence rates for the 
conformal and adapted grids.

\subsection{Grid quality}
Now we construct our elliptic grids and compare them to 
the near-conformal grid suggested by Reference~\cite{Ribeiro2010} as a reference grid for existing flux aligned grids.
Note that in this section by orthogonal grid we always mean the adapted orthogonal grid discussed in Section~\ref{sec:adaption}. 
\begin{figure}[htbp]
\centering
\subfloat[orthogonal $\psi_0=-20$]{\includegraphics[trim = 0px 0px 0px 0px, clip, scale=0.3]{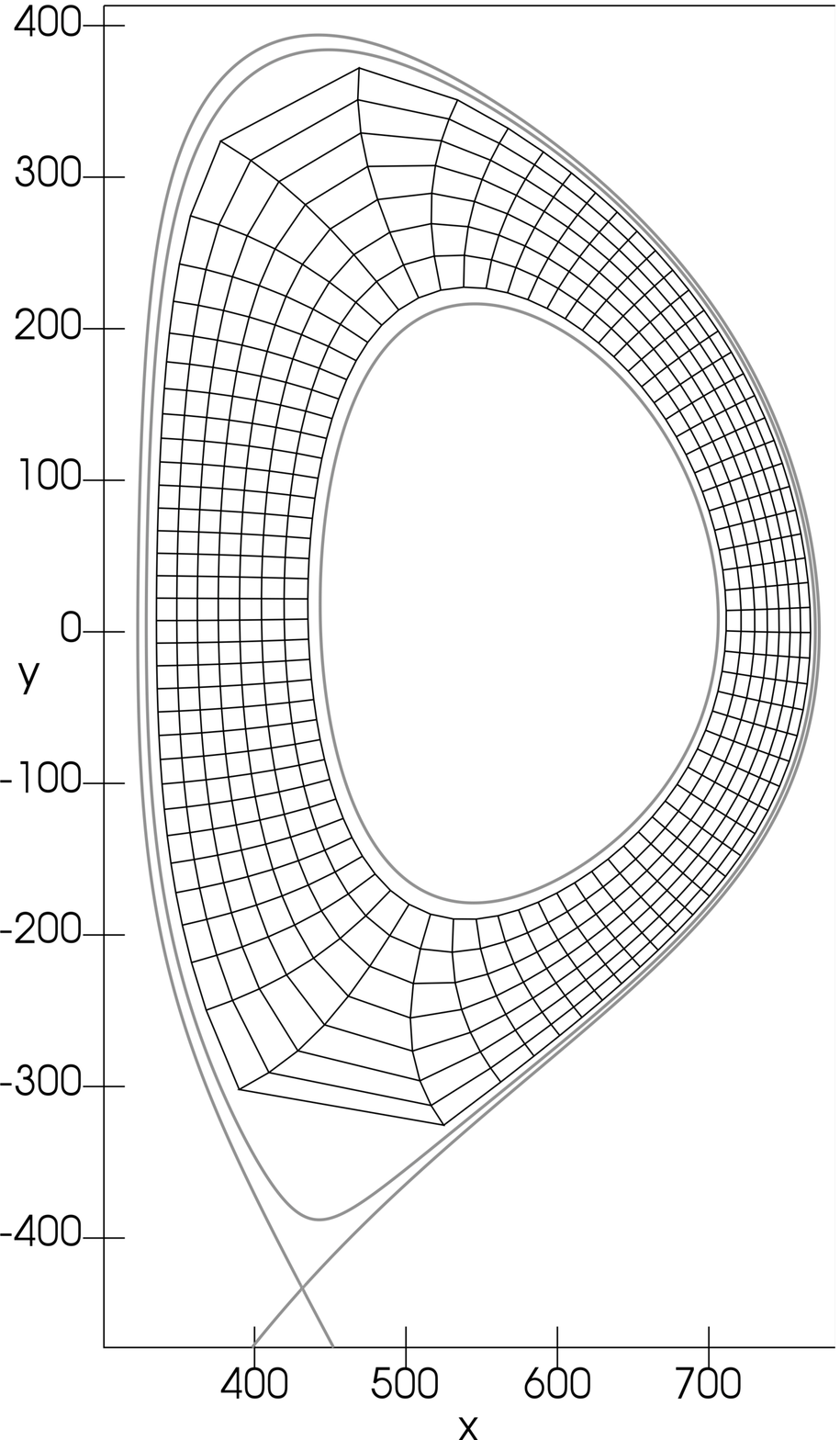}\label{fig:three_gridsa}}
\subfloat[orthogonal $\psi_0=-1$]{\includegraphics[trim = 0px 0px 0px 0px, clip, scale=0.3]{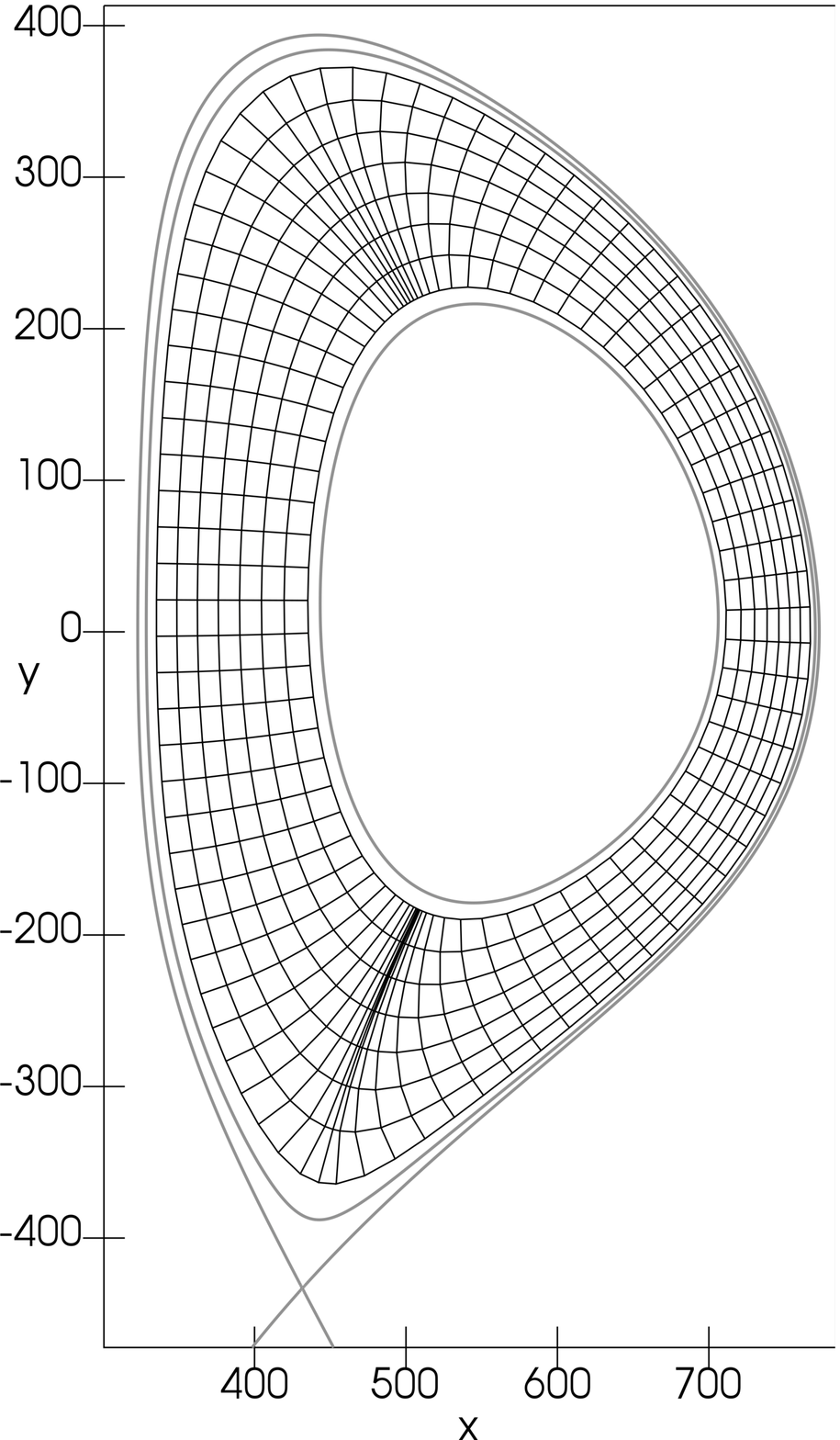}\label{fig:three_gridsb}}
\subfloat[near conformal]{\includegraphics[trim = 0px 0px 0px 0px, clip, scale=0.3]{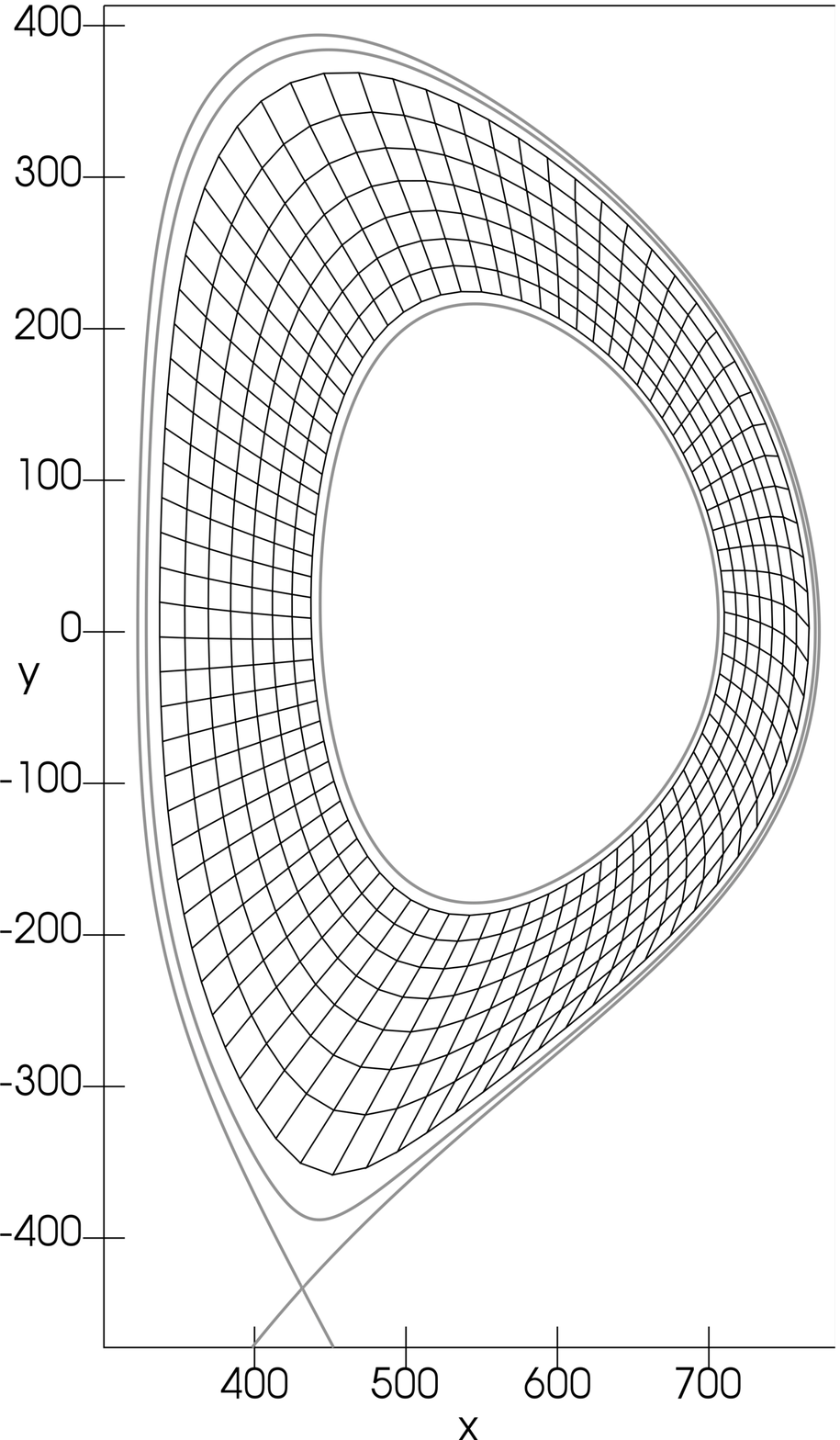}\label{fig:three_gridsc}}
\caption{
  Orthogonal grid with $\psi_0=-20$, $\psi_1=-1$ (a), orthogonal grid with $\psi_0=-1$, $\psi_1=-20$  (b) both with $w=|\nabla\psi|$ and near conformal grid with $\psi_0=-20$, $\psi_1=-1$ (c) with equal number of grid points 
  $P=1$, $N_\zeta=8$, $N_\eta=80$. 
Note that the nodes represent cell centers and not the actual cell boundaries. 
The grey lines denote the actual grid boundaries at $\psi=-20$ and $\psi=-1$. For 
the sake of orientation we also plot the line $\psi=0$ to indicate the X-point. 
}
\label{fig:three_grids}
\end{figure}
In Fig.~\ref{fig:three_grids} we plot the orthogonal grid described in Section \ref{sec:geometry} with the function $h$ chosen once on the
inner boundary~\ref{fig:three_gridsa} and once on the outer boundary~\ref{fig:three_gridsb} 
(in both cases adaption is used), and the near conformal grid proposed by Reference~\cite{Ribeiro2010} in~\ref{fig:three_gridsc}.
In the orthogonal grids, we have an equidistant discretization at the inner respectively outer boundary. Each of these grids is aligned with the magnetic flux function
$\psi$. 
The orthogonal grids \ref{fig:three_gridsa} and \ref{fig:three_gridsb} distribute equally spaced cells evenly in the region around $Z=0$,
where the contour lines are almost straight. 
The cell size in the radial direction decreases as we go from the inner to the outer boundary. 
However, when the curvature of the
contour lines increases (this is the case in the region around $R=500$), the cells are either  
 prolongated or compressed in the direction of the angle-like coordinate $\eta$. 
 This effect results in cells with very large aspect ratio of the cells in the
 orthogonal grids.

 The near-conformal grid~\ref{fig:three_gridsc} is near-conformal in the sense that the conformal deformation $g^{\zeta\zeta}/g^{\eta\eta}$ is small. Note, however, that this grid is
completely aligned to the function $\psi$. 
This means that the aspect ratio of the cells remains
as constant as possible. 
The downside of this grid is that it is clearly non-orthogonal especially close to the boundary and thus Neumann boundary conditions are difficult to implement.   

\begin{figure}[htbp]
\centering
\subfloat[true conformal grid]{\includegraphics[trim = 0px 0px 0px 0px, clip, scale=0.3]{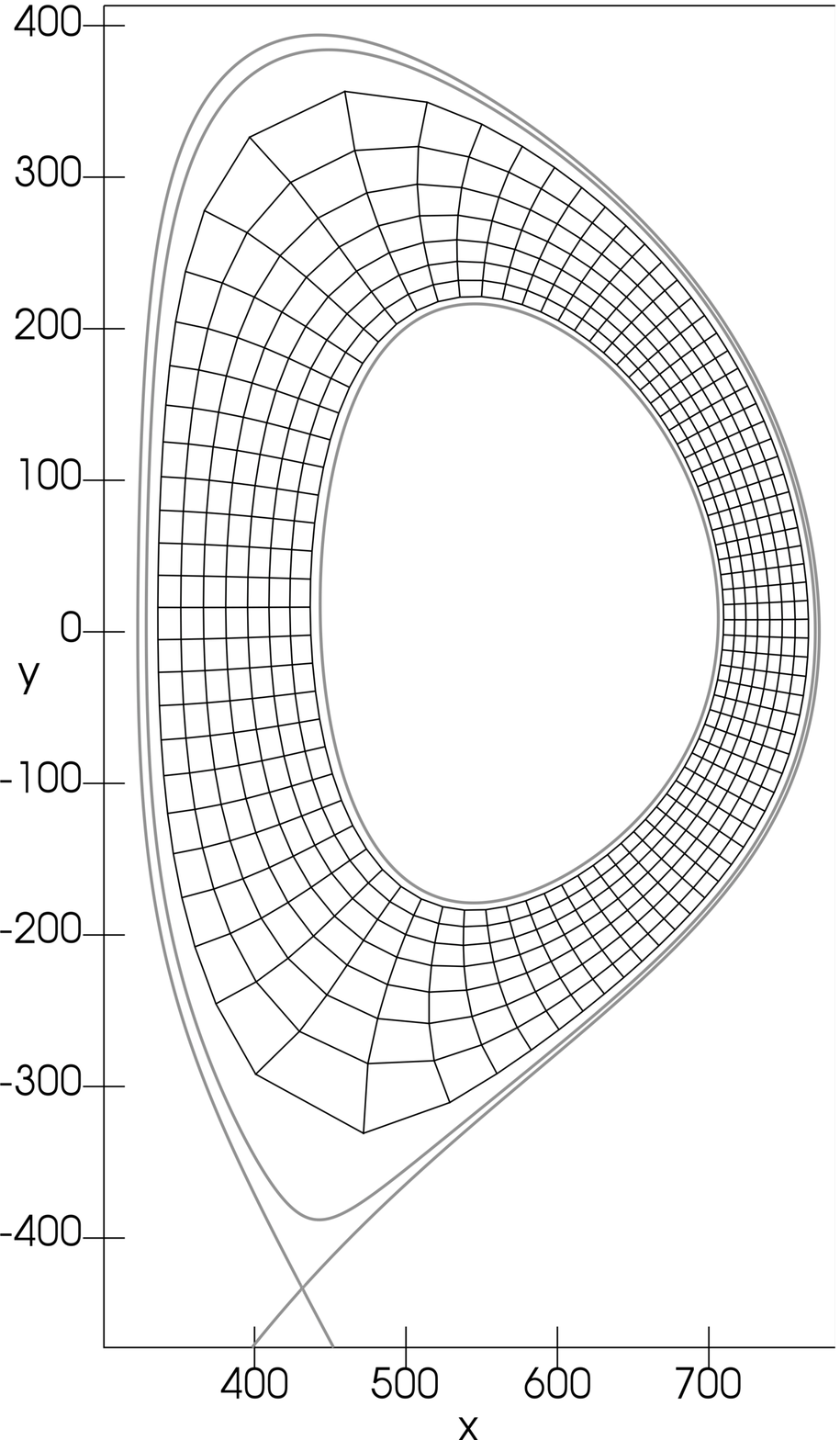}\label{fig:modified_gridsa}}
\subfloat[grid adaption]{\includegraphics[trim = 0px 0px 0px 0px, clip, scale=0.3]{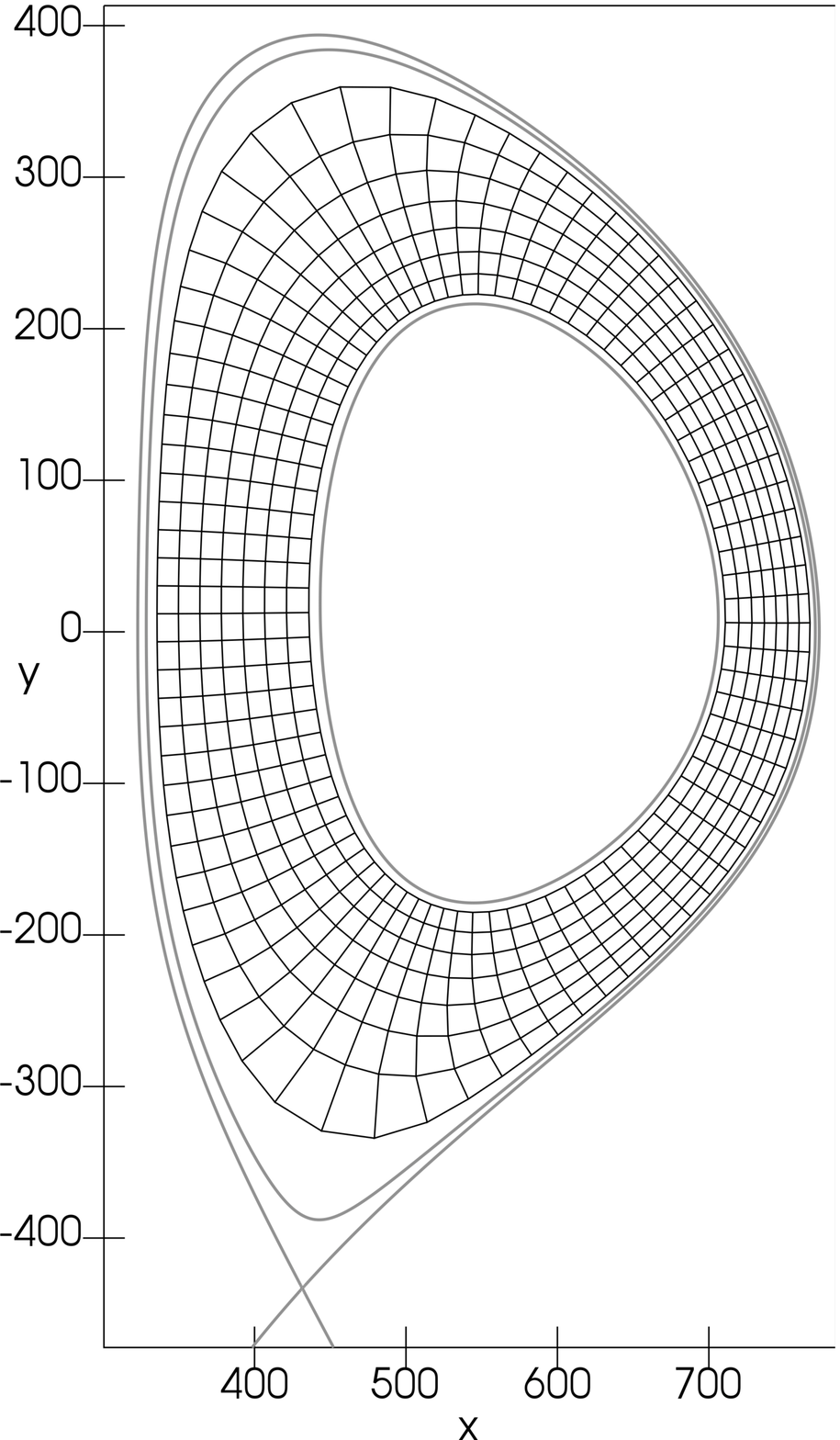}\label{fig:modified_gridsb}}
\subfloat[monitor metric]{\includegraphics[trim = 0px 0px 0px 0px, clip, scale=0.3]{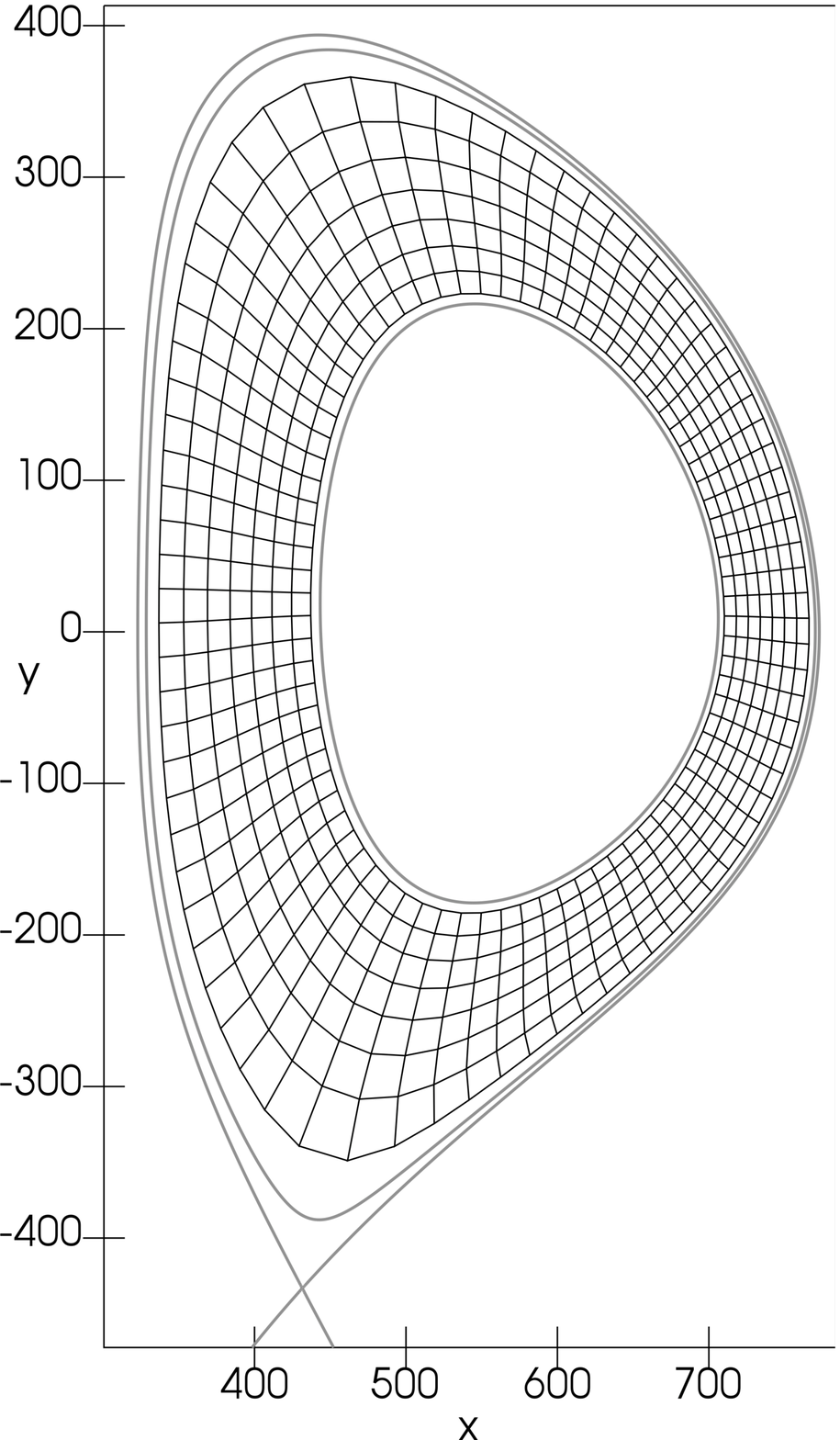}\label{fig:modified_gridsc}}
\caption{
  Truly conformal grid (a), adapted grid with $w=|\nabla\psi|$ (b) and monitor grid with $k=0.1$ and $\eps=0.001$ (c) with equal number of grid points 
  $P=1$, $N_\zeta=8$, $N_\eta=80$. We choose $\psi_0 = -20$ and $\psi_1 = -1$. 
Note that the nodes represent cell centers and not the actual cell boundaries. 
The grey lines denote the actual grid boundaries at $\psi=-20$ and $\psi=-1$. For 
the sake of orientation we also plot the line $\psi=0$ to indicate the X-point. 
}
\label{fig:modified_grids}
\end{figure}

In Fig.~\ref{fig:modified_grids} we plot the results from 
our elliptic grid generation processes. We show a 
true conformal grid~\ref{fig:modified_gridsa}, an elliptic grid with adaption~\ref{fig:modified_gridsb} and 
one constructed with a monitor metric~\ref{fig:modified_gridsc}. The last grid is clearly 
non-orthogonal while the conformal and adapted grids are orthogonal in the whole domain. 
All three grids are orthogonal at the boundary (which facilitates the implementation of Neumann boundary conditions). 
The conformal grid~\ref{fig:modified_gridsa} suffers from large cells in the upper 
and lower regions of the domain. This is cured 
by either performing grid adaption~\ref{fig:modified_gridsb} or by employing a monitor metric~\ref{fig:modified_gridsc}.
Both result in cells of similar size. The difference
between these grids is the degree of flux alignment, which is clearly superior
in the monitor grid in the region close to the X-point (Reference~\cite{Glasser2006} provides a more quantitative investigation of the degree of alignment).

We assess the quality of our grids using an analytic solution of the following elliptic equation
\begin{equation}
    \nabla\cdot\left( \chi\nabla f \right) = \rho, \label{eq:analytic}
\end{equation}
where we either choose a flux aligned solution
\begin{subequations}
\begin{align}
f_{\text{ana}}(x,y)=& 0.1(\psi-\psi_0)(\psi-2\psi_1+\psi_0) \\
\chi(x,y) =& \frac{x_0}{x} \sqrt{1 +(\nabla \psi)^2}
  \label{}
\end{align}
\label{eq:analytic_field}
\end{subequations}
or a localized solution
\begin{subequations}
\begin{align}
f_{\text{ana}}(x,y)=& 
\begin{cases}
0 \text{ if } \left( \frac{x-x_b}{\sigma} \right)^2+\left( \frac{y-y_b}{\sigma} \right)^2 > 1 \\
\exp \left[ 1+\left( \left( \frac{x-x_b}{\sigma} \right)^2+\left( \frac{y-y_b}{\sigma} \right)^2 - 1 \right)^{-1} \right] \text{ else }
\end{cases} \\
\chi(x,y) =& \frac{x_0}{x} \sqrt{1 +(\nabla \psi)^2}(1+0.5\sin(\theta(x,y))).
  \label{}
\end{align}
\label{eq:analytic_local}
\end{subequations}
The corresponding parameters are as follows: $x_0 = 550$, $\psi_0 = -20$ and $\psi_1=-1$, $x_b=440$, $y_b=-220$ and $\sigma=40$. 
The position of the blob corresponds to the lower left corner of the domain. In this region the grids differ the most.
Also note that the flux aligned solution~\eqref{eq:analytic_field} 
has a Neumann boundary condition at $\psi=\psi_1$ and a Dirichlet
boundary at $\psi=\psi_0$. For the localized solution~\eqref{eq:analytic_local} we choose
Dirichlet conditions at both boundaries.
In order to numerically solve Eq.~\eqref{eq:analytic} we first have to transform the metric
and the involved functions to the coordinate system that we use (cf. Reference~\cite{Einkemmer2014, Held2016} on how the derivatives are discretized).
It might seem that the conformal and orthogonal grids have a 
performance advantage over the curvilinear grids because the non-diagonal elements
of the metric vanish. 
However, the number of necessary matrix-vector multiplications is four in all cases and only the number of vector-vector multiplications is different. Consequently, we only observe minor differences in performance. 

The relative errors are computed in the $L2$ norm via
\begin{align}
  \eps = \left(\frac{\int_{u_0}^{u_1}\d u \int_0^{2\pi}\d v \sqrt{g} (f_{\text{num}} - f_{\text{ana}})^2 }
  {\int_{u_0}^{u_1}\d u \int_0^{2\pi}\d v \sqrt{g} f_{\text{ana}}^2 }\right)^{1/2}
  \label{eq:error_norm}
\end{align}
where $\sqrt{g}\d u \d v$ is the correct volume form in the $u,v$ coordinate system. 
The ratio of $N_v/N_u$ is chosen such that the aspect 
ratio of the resulting cells is approximately unity.

\begin{table}[htbp]
  \small
  \begin{tabular}{|c|c|c|c|c|c|c|c|c|c|}
    \hline
   $N_u$ & $N_v$   & \multicolumn{2}{c|}{Conformal} &  \multicolumn{2}{c|}{Adapted} & \multicolumn{2}{c|}{Monitor} & \multicolumn{2}{c|}{Orthogonal}  \\ \hline
 \multicolumn{10}{|c|}{P=3} \\ \hline
2 & 20 & 3.10E-02 & - & 1.37E-02 & - & 8.15E-03 & - & 7.34E-03 & - \\ \hline
4 & 40 & 9.58E-03 & 1.69 & 4.80E-03 & 1.51 & 1.86E-03 & 2.13 & 4.67E-04 & 3.97 \\ \hline
8 & 80 & 1.32E-03 & 2.86 & 1.17E-03 & 2.03 & 2.70E-04 & 2.78 & 4.28E-05 & 3.45 \\ \hline
16 & 160 & 7.76E-04 & 0.76 & 2.74E-04 & 2.10 & 4.42E-05 & 2.61 & 4.26E-06 & 3.33 \\ \hline
32 & 320 & 5.17E-05 & 3.91 & 5.02E-05 & 2.45 & 6.14E-06 & 2.85 & 3.58E-07 & 3.57 \\ \hline
64 & 640 & 4.57E-06 & 3.50 & 7.41E-06 & 2.76 & 7.85E-07 & 2.97 & 3.23E-08 & 3.47 \\ \hline
 \multicolumn{10}{|c|}{P=4}  \\ \hline
2 & 20 & 5.30E-03 & - & 8.72E-03 & - & 2.72E-03 & - & 3.01E-04 & - \\ \hline
4 & 40 & 1.48E-03 & 1.84 & 1.61E-03 & 2.43 & 2.21E-04 & 3.62 & 1.32E-05 & 4.51 \\ \hline
8 & 80 & 9.30E-04 & 0.67 & 3.87E-04 & 2.06 & 5.47E-05 & 2.01 & 1.10E-06 & 3.59 \\ \hline
16 & 160 & 1.31E-04 & 2.83 & 6.96E-05 & 2.48 & 7.57E-06 & 2.85 & 4.14E-08 & 4.73 \\ \hline
32 & 320 & 4.57E-06 & 4.84 & 1.00E-05 & 2.80 & 7.83E-07 & 3.27 & 2.01E-09 & 4.37 \\ \hline
64 & 640 & 5.66E-07 & 3.01 & 1.07E-06 & 3.22 & 7.93E-08 & 3.30 & 9.61E-11 & 4.39 \\ \hline
\end{tabular}
    \caption{ Convergence table for Eq.~\eqref{eq:analytic} with the flux aligned 
    solution given in Eq.~\eqref{eq:analytic_field} on the different 
    grids proposed.
  We show the relative error \eqref{eq:error_norm} and the corresponding order 
  for various numbers of polynomial coefficients $P$ in the dG discretization. 
  The orthogonal grid uses $\psi_0=-20$ as the first line.
}
\label{tab:fluxaligned}
\end{table}

Let us point out first that if we choose a function $f_{\text{ana}}=f_{\text{ana}}(\psi)$ the 
flux aligned grids produce errors that can be orders of magnitude lower
than those of the non-aligned grids. This is observed in Table~\ref{tab:fluxaligned}.
The orthogonal grid has much better errors than the elliptic grids. 
However, while the conformal and adapted grids have similar errors the monitor grid shows smaller errors for all resolutions. 
This is due to the fact that the degree of alignment is higher in case of the monitor grid (see, for example, \cite{Glasser2006} for a discussion).
If due to some physical reasons the solution is expected to be flux aligned,
the aligned grids are thus preferable. 
From the dG method we expect convergence of order $P$~\cite{Arnold2001}.
This is achieved only in ideal cases, however. We observe more irregular behaviour. The orthogonal grid shows approximate orders of $3.5$ and $4.5$, while 
the remaining grids show orders mostly between $2.0$ and $3.5$. 
Since the discretization method itself is equal in all four cases this 
irregularity is most likely due to the different metric elements in the
four coordinate systems.

In actual turbulence simulations we expect localized structures and eddies
as exemplified by the localized solution given in Eq.~\eqref{eq:analytic_local}.
\begin{table}[htbp]
  \small
  \begin{tabular}{|c|c|c|c|c|c|c|c|c|c|}
    \hline
   $N_u$ & $N_v$   & \multicolumn{2}{c|}{Orthogonal } &  \multicolumn{2}{c|}{Conformal} & \multicolumn{2}{c|}{Near Conformal} & \multicolumn{2}{c|}{Monitor}  \\ \hline
 \multicolumn{10}{|c|}{P=3} \\ \hline
4 & 40 & 1.24E+00 & - & 1.07E+00 & - & 6.67E+00 & - & 2.37E+00 & - \\ \hline
8 & 80 & 3.38E-01 & 1.87 & 7.39E-01 & 0.54 & 5.35E-01 & 3.64 & 6.53E-01 & 1.86 \\ \hline
16 & 160 & 2.57E-01 & 0.39 & 3.20E-01 & 1.21 & 2.92E-01 & 0.87 & 1.15E-01 & 2.51 \\ \hline
32 & 320 & 3.68E-02 & 2.80 & 6.98E-02 & 2.20 & 1.76E-02 & 4.06 & 1.38E-02 & 3.06 \\ \hline
64 & 640 & 3.15E-03 & 3.55 & 6.70E-04 & 6.70 & 1.58E-03 & 3.47 & 3.78E-04 & 5.18 \\ \hline
 \multicolumn{10}{|c|}{P=4}  \\ \hline
4 & 40 & 2.39E+00 & - & 1.49E+00 & - & 1.84E+00 & - & 1.19E+00 & - \\ \hline
8 & 80 & 6.41E-01 & 1.90 & 8.11E-01 & 0.88 & 7.59E-01 & 1.28 & 4.09E-01 & 1.54 \\ \hline
16 & 160 & 1.18E-01 & 2.44 & 1.04E-01 & 2.96 & 9.81E-02 & 2.95 & 5.56E-02 & 2.88 \\ \hline
32 & 320 & 1.56E-02 & 2.92 & 2.14E-02 & 2.28 & 7.66E-03 & 3.68 & 1.93E-03 & 4.85 \\ \hline
64 & 640 & 6.78E-04 & 4.52 & 3.66E-04 & 5.87 & 2.62E-04 & 4.87 & 5.15E-05 & 5.23 \\ \hline
\end{tabular}
    \caption{ Convergence table for the solution of Eq.~\eqref{eq:analytic} with the localized solution given by Eq.~\eqref{eq:analytic_local} on the different grids proposed.
  We show the relative error \eqref{eq:error_norm} and the corresponding order 
  for two different numbers of polynomial coefficients $P$ in the dG discretization. The orthogonal grid uses $\psi_0=-20$ as the first line.
}
\label{tab:elliptic}
\end{table}

In Table \ref{tab:elliptic} we show the resulting convergence for a forward discretization on the orthogonal grid, the truly conformal grid, the near conformal grid, 
and the elliptic grid with monitor metric. 

Note that we again observe irregular convergence behaviour for all grids. Nonetheless we ordered the 
grids from mostly high to lower errors from left to right. 
The grid constructed with the monitor metric has the lowest error. The worst grids are the orthogonal and the true conformal
grids, which have the largest cell sizes at the location,
where the blob was placed. 
The errors from the adapted grid are 
approximately between the near conformal grid and the grid 
with monitor metric. 

The performance of these grids can be more clearly understood by comparing their maximal and minimal cell sizes for a fixed number of cells. Recall that the arc length of a curve is the integral over the norm 
of the tangent vector of the curve. In the case of the coordinate
lines the tangent vectors are just the basis vectors $\partial_u$ and $\partial_v$.
If we assume that the computational space is discretized in equidistant
cells with cell size $h_u$ and $h_v$ respectively,
we define
\begin{subequations}
\begin{align}
  l_u := \sqrt{g_{uu}} h_u = \sqrt{g}\sqrt{ g^{vv}} h_u \\
  l_v := \sqrt{g_{vv}} h_v = \sqrt{g}\sqrt{ g^{uu}} h_v
\end{align}
  \label{eq:cellsize}
\end{subequations}
Furthermore we define the ratio of maximal to minimal cell size by
\begin{align}
  a_u := \frac{\max l_u}{\min l_u} \quad
  a_v := \frac{\max l_v}{\min l_v}
  \label{eq:ratios}
\end{align}
\begin{table}[htbp]
\begin{tabular}{|l|c|c|c|c||c|c|}
\hline
 & {$\max l_u$} & {$\max l_v$} & {$\min l_u$} & $\min l_v$ & $a_u$ & $a_v$\\ \hline
Orthogonal ($\psi_0=-1$) & 14.23 & 6.03 & 1.53 & 0.06 & 9.28 & 94.56 \\ \hline
Orthogonal ($\psi_0=-20$) & 9.50 & 154.22 & 1.53 & 3.28 & 6.20 & 47.09 \\ \hline
Conformal & 57.02 & 79.51 & 1.75 & 2.44 & 32.53 & 32.53 \\ \hline
Near Conformal & 21.99 & 31.07 & 1.62 & 2.28 & 13.59 & 13.62 \\ \hline
Adapted & 58.91 & 20.05 & 1.68 & 2.04 & 34.98 & 9.82 \\ \hline
Monitor & 33.14 & 15.42 & 1.96 & 3.04 & 16.91 & 5.07 \\ \hline
\end{tabular}
\caption{Comparison of minimal and maximal cell sizes as defined in Eq.~\eqref{eq:cellsize} and the ratios defined in Eq.~\eqref{eq:ratios} of the various grids with equal number of gridpoints $N_u=32$, $N_v=320$. The grids
are ordered according to the ratio of maximum to minimum cell size in $v$ from top to bottom. }
\label{tab:gridsize}
\end{table}

In Table~\ref{tab:gridsize} we compute the maximum and minimum cell
sizes of the various grids. 
For completeness we also computed the values for the orthogonal grid with the 
first line being the $\psi_0=-1$ line, which is plotted in Fig.~\ref{fig:three_grids}.
We expect that if the maximal cell size is small, then the error committed by this grid is also small. 
In fact, this explains the results that were obtained in Table \ref{tab:elliptic}. 
That is, the monitor grid and the near conformal grid are better in terms of accuracy compared to the conformal and orthogonal grids. 
Note that the orthogonal grid with $\psi_0=-20$ as the fist line has quite long cells in the $v$ direction but smaller cells than all the other grids in the $u$ direction.

The accuracy obtained is an important consideration. 
However, it has to be seen in the context that this elliptic equation will eventually be coupled to an evolution equation. 
In this case the minimal grid size is an important characteristics as it determines the CFL condition of explicit solvers and the 
condition number of the discretization matrix in implicit schemes
and thus the computational efficiency of the grid used. 
Consider the minimal grid size in the $v$ direction of the orthogonal grid with $\psi_0=-1$ as the first line. It is $50$ times  smaller than the cell in the monitor grid.
Thus, this grid is clearly not efficient for a practical numerical simulation. 

The values of $a_u$ and $a_v$ are a measure of how well a grid performs. 
Small values indicate small errors in the elliptic equation and a large
minimal cell size, which is advantageous for advection schemes. 
Overall we thus find that the monitor grid gives the best accuracy as well as the least stringent CFL condition among the grids we have considered here.
This is expressed by a minimal value of $a_v$ and the $a_u$ value second only to
the near conformal grid.

\section{Conclusion} \label{sec:conclusion}
All in all, we show that we are able to construct various elliptic grids by the method
of streamline integration combined with the solution of a suitably chosen elliptic equation.
Compared to the TTM and related methods~\cite{Thompson1977, Glasser2006, Vaseva2009} we significantly simplify the construction and
accuracy of the method. Compared to the flux-aligned grids suggested in the literature
on magnetic fusion problems~\cite{Ribeiro2010} we improve the distribution of cells across the 
domain and enable boundary orthogonality.

In our example implementation we discretize the approach laid out in Section~\ref{sec:geometry} with high order Runge Kutta and dG methods. 
Our grids are suitable but not limited to the discretization of the edge region of magnetically
confined plasmas as proved in Section~\ref{sec:numerics}. We explicitly show convergence of two solutions for a general elliptic 
equation that appears in typical physical models of plasma turbulence. 
Our orthogonal grids are flux aligned but show the largest 
variation in cell size of all the grids we investigated. 
The conformal grid to a lesser extent suffers from the same problem.
In conclusion we find that the adaptive grid and even more the monitor metric approach yield the best grids in terms of small error and homogeneous cell distribution.   

In the future we intend to investigate the possibility to include the X-point into
the domain of interest. The difficulty encountered with our algorithm at
an X-point or O-point is the vanishing gradient $\nabla\psi$. The Jacobian of 
any flux aligned coordinate system thus has a singularity at this point
and the convergence of elliptic equations significantly deteriorates or even 
vanishes. Numerical methods that can handle singularities might solve the problem.
Furthermore, a generalization of our algorithm to 
three dimensions might be feasible.

\appendix
\section{Conformal and field aligned coordinates \label{app:conf-field}}

The magnetic flux $\psi$ does in general not satisfy the Laplace equation and thus is not able to serve as a conformal coordinate directly. It is therefore unclear whether it is possible to find a field aligned conformal coordinate system. Thus, we ask if it is possible to find a function $u$ such that $u(\psi)$ satisfies the Laplace equation
\[ \Delta u(\psi) = 0 \]
which is equivalent to
\[ (\nabla\psi)^{2}u^{\prime\prime}(\psi)+(\Delta\psi)u^{\prime}(\psi) = 0 \]
Now we can rewrite this equation as
\[ u^{\prime\prime}(\psi)=-\frac{\Delta\psi}{(\nabla\psi)^{2}}u^{\prime}(\psi). \]
Since the left-hand side only depends on $\psi$, the same must hold true for the right-hand side. This implies that
\begin{equation} \frac{\Delta\psi}{(\nabla\psi)^{2}}=g(\psi), \label{eq:conform-condition} \end{equation}
where $g$ is an arbitrary function of $\psi$. Note that strictly speaking it is also possible to have $u^{\prime\prime}=u^{\prime}=0$ which implies that $u(\psi)=\text{const}$. This is certainly a solution but can be discarded for the purpose of constructing a coordinate system.

Equation \eqref{eq:conform-condition} gives a condition on $\psi$ that, if satisfied, allows us to construct a conformal field aligned coordinate system. Unfortunately, this is not a property that holds for the solutions of the Grad--Shafranov equation in general:
\begin{align}
  R\frac{\partial}{\partial R}\left( \frac{1}{R}\frac{\partial \psi}{\partial R} \right) + \frac{\partial \psi}{\partial Z^2} = -\mu_0 R \frac{\d p(\psi)}{\d \psi} -F(\psi)\frac{\d F(\psi)}{\d \psi}
  \label{eq:GS}
\end{align}
For example, $\psi(R,Z)=R^2Z^2$ is an equilibrium solution such that
\[ \frac{\Delta\psi}{(\nabla\psi)^{2}} = \frac{1}{2\psi}\left(1+\frac{1}{1+R^{2}/Z^{2}}\right) \]
which clearly does not satisfy condition~\eqref{eq:conform-condition}.

Note, however, that for specific equilibria the condition given by equation \eqref{eq:conform-condition} can be satisfied. For example, $\psi(R,Z)=R^4$ satisfies the Grad--Shafranov equation and we have
\[ \frac{\Delta\psi}{(\nabla\psi)^{2}}=\frac{1}{\psi}. \]
Both this and the equilibrium considered above are particular cases of the class of Solov\'ev equilibria.
Although the example is uninteresting for the purpose of turbulence simulations it is noteworthy 
that solutions that satisfy both~\eqref{eq:conform-condition} and \eqref{eq:GS} do exist.

\section{Coefficients \label{app:coeff-psi}}
For ease of reproduction we print the coefficients of the solov\'ev
equilibrium solution for $\psi$ in Reference~\cite{Cerfon2010}
\begin{align*}
    A &= 0.0,\\
    c_{1..12} &=[0.07350114445500399706,
           -0.08662417436317227513,\\
           &\quad-0.14639315434011026207,
           -0.07631237100536276213,\\
            &\quad0.09031790113794227394,
           -0.09157541239018724584,\\
           &\quad-0.003892282979837564482,
         0.04271891225076417603, \\
            &\quad0.22755456460027913117,
          -0.13047241360177695448, \\
           &\quad-0.03006974108476955225,
            0.004212671892103931173 ],\\
    R_0   &= 547.891714877869, \\
    a/R_0 &= 0.41071428571428575, \\
    \kappa     &= 1.75,\\
    \delta     &= 0.47,\\
\end{align*}
These coefficients are also contained in the dataset for this paper~\cite{FELTORv3.1}.

\section*{Acknowledgements} 	
This work was supported by the Austrian Science Fund (FWF) W1227-N16 and Y398. 
This work has been carried out within the framework of the EUROfusion Consortium and has received funding from the Euratom research and training programme 2014‐2018 under grant agreement No 633053. The views and opinions expressed herein do not necessarily reflect those of the European Commission.

\bibliography{conformal}

\begin{thebibliography}{27}
\expandafter\ifx\csname natexlab\endcsname\relax\def\natexlab#1{#1}\fi
\providecommand{\bibinfo}[2]{#2}
\ifx\xfnm\relax \def\xfnm[#1]{\unskip,\space#1}\fi
\bibitem[{Thompson et~al.(1997)Thompson, Warsi, and Mastin}]{Thompson}
\bibinfo{author}{J.~E. Thompson}, \bibinfo{author}{Z.~Warsi},
  \bibinfo{author}{C.~W. Mastin}, \bibinfo{title}{Numerical Grid Generation},
  \bibinfo{publisher}{Joe F. Thompson}, \bibinfo{year}{1997}.
\bibitem[{Weatherill et~al.(1998)Weatherill, Soni, and Thompson}]{Grid}
\bibinfo{author}{N.~Weatherill}, \bibinfo{author}{B.~Soni},
  \bibinfo{author}{J.~Thompson}, \bibinfo{title}{Handbook of Grid Generation},
  \bibinfo{publisher}{Informa {UK} Limited}, \bibinfo{year}{1998}.
\bibitem[{Liseikin(2007)}]{Liseikin}
\bibinfo{author}{V.~D. Liseikin}, \bibinfo{title}{A Computational Differential
  Geometry Approach to Grid Generation}, \bibinfo{publisher}{Springer-Verlag},
  \bibinfo{edition}{second} edition, \bibinfo{year}{2007}.
\bibitem[{Thompson et~al.(1977)Thompson, Thames, and Mastin}]{Thompson1977}
\bibinfo{author}{J.~F. Thompson}, \bibinfo{author}{F.~C. Thames},
  \bibinfo{author}{C.~W. Mastin},
\newblock \bibinfo{title}{Tomcat - code for numerical generation of
  boundary-fitted curvilinear coordinate systems on fields containing any
  number of arbitrary 2-dimensional bodies},
\newblock \bibinfo{journal}{J. Comput. Phys.} \bibinfo{volume}{24}
  (\bibinfo{year}{1977}) \bibinfo{pages}{274--302}.
\bibitem[{Papamichael and Stylianopoulos(2010)}]{Papamichael}
\bibinfo{author}{N.~Papamichael}, \bibinfo{author}{N.~Stylianopoulos},
  \bibinfo{title}{Numerical Conformal Mapping}, \bibinfo{publisher}{World
  Scientific Publishing}, \bibinfo{year}{2010}.
\bibitem[{Glasser et~al.(2006)Glasser, Liseikin, Vaseva, and
  Likhanova}]{Glasser2006}
\bibinfo{author}{A.~H. Glasser}, \bibinfo{author}{V.~D. Liseikin},
  \bibinfo{author}{I.~A. Vaseva}, \bibinfo{author}{Y.~V. Likhanova},
\newblock \bibinfo{title}{Some computational aspects on generating numerical
  grids},
\newblock \bibinfo{journal}{Russ. J. Numer. Anal. Math. Modelling}
  \bibinfo{volume}{21} (\bibinfo{year}{2006}) \bibinfo{pages}{481--505}.
\bibitem[{Vaseva et~al.(2009)Vaseva, Liseikin, Likhanova, and
  Morokov}]{Vaseva2009}
\bibinfo{author}{I.~A. Vaseva}, \bibinfo{author}{V.~D. Liseikin},
  \bibinfo{author}{Y.~V. Likhanova}, \bibinfo{author}{Y.~N. Morokov},
\newblock \bibinfo{title}{An elliptic method for construction of adaptive
  spatial grids},
\newblock \bibinfo{journal}{Russ. J. Numer. Anal. Math. Modelling}
  \bibinfo{volume}{24} (\bibinfo{year}{2009}) \bibinfo{pages}{65--78}.
\bibitem[{Wesson(2011)}]{Wesson}
\bibinfo{author}{J.~Wesson}, \bibinfo{title}{Tokamaks},
  \bibinfo{publisher}{Oxford University Press}, \bibinfo{edition}{4th} edition,
  \bibinfo{year}{2011}.
\bibitem[{Hamada(1962)}]{hamada62}
\bibinfo{author}{S.~Hamada},
\newblock \bibinfo{title}{Hydromagnetic equilibria and their proper
  coordinates},
\newblock \bibinfo{journal}{Nucl. Fusion} \bibinfo{volume}{2}
  (\bibinfo{year}{1962}).
\bibitem[{Boozer(1980)}]{boozer80}
\bibinfo{author}{A.~H. Boozer},
\newblock \bibinfo{title}{{Guiding center drift equations}},
\newblock \bibinfo{journal}{Phys. Fluids}  (\bibinfo{year}{1980})
  \bibinfo{pages}{904}.
\bibitem[{Boozer(1981)}]{boozer81}
\bibinfo{author}{A.~H. Boozer},
\newblock \bibinfo{title}{{Plasma equilibrium with rational magnetic
  surfaces}},
\newblock \bibinfo{journal}{Phys. Fluids}  (\bibinfo{year}{1981})
  \bibinfo{pages}{1999}.
\bibitem[{Grimm et~al.(1983)Grimm, Dewar, and J.}]{Grimm1983}
\bibinfo{author}{R.~Grimm}, \bibinfo{author}{R.~Dewar},
  \bibinfo{author}{M.~J.},
\newblock \bibinfo{title}{{Ideal MHD stability calculations in axisymmetric
  toroidal coordinate systems}},
\newblock \bibinfo{journal}{J. Comput. Phys.} \bibinfo{volume}{48}
  (\bibinfo{year}{1983}) \bibinfo{pages}{94--117}.
\bibitem[{Cheng(1992)}]{cheng92}
\bibinfo{author}{C.~Z. Cheng},
\newblock \bibinfo{title}{{Kinetic extensions of magnetohydrodynamics for
  axisymmetric toroidal plasmas }},
\newblock \bibinfo{journal}{Physics Reports} \bibinfo{volume}{211}
  (\bibinfo{year}{1992}) \bibinfo{pages}{1--51}.
\bibitem[{Park et~al.(2008)Park, Boozer, and Menard}]{park08}
\bibinfo{author}{J.-k. Park}, \bibinfo{author}{A.~H. Boozer},
  \bibinfo{author}{J.~E. Menard},
\newblock \bibinfo{title}{{Spectral asymmetry due to magnetic coordinates}},
\newblock \bibinfo{journal}{Phys. Plasmas}  (\bibinfo{year}{2008})
  \bibinfo{pages}{064501}.
\bibitem[{Czarny and Huysmans(2008)}]{Czarny2008}
\bibinfo{author}{O.~Czarny}, \bibinfo{author}{G.~Huysmans},
\newblock \bibinfo{title}{Bezier surfaces and finite elements for mhd
  simulations},
\newblock \bibinfo{journal}{J. Comput. Phys.} \bibinfo{volume}{227}
  (\bibinfo{year}{2008}) \bibinfo{pages}{7423--7445}.
\bibitem[{Ribeiro and Scott(2010)}]{Ribeiro2010}
\bibinfo{author}{T.~T. Ribeiro}, \bibinfo{author}{B.~D. Scott},
\newblock \bibinfo{title}{Conformal tokamak geometry for turbulence
  computations},
\newblock \bibinfo{journal}{IEEE Transactions on Plasma Science}
  \bibinfo{volume}{38} (\bibinfo{year}{2010}) \bibinfo{pages}{2159--2168}.
\bibitem[{Oliver et~al.(2012)Oliver, Montero, Montenegro, Rodr\'{\i}guez,
  Escobar, and Perez-Foguet}]{Oliver2012}
\bibinfo{author}{A.~Oliver}, \bibinfo{author}{G.~Montero},
  \bibinfo{author}{R.~Montenegro}, \bibinfo{author}{E.~Rodr\'{\i}guez},
  \bibinfo{author}{J.~M. Escobar}, \bibinfo{author}{A.~Perez-Foguet},
\newblock \bibinfo{title}{Finite element simulation of a local scale air
  quality model over complex terrain},
\newblock \bibinfo{journal}{Adv. Sci. Res.} \bibinfo{volume}{8}
  (\bibinfo{year}{2012}) \bibinfo{pages}{105--113}.
\bibitem[{D'haeseleer et~al.(1991)D'haeseleer, Hitchon, Callen, and
  Shohet}]{haeseleer}
\bibinfo{author}{W.~D'haeseleer}, \bibinfo{author}{W.~Hitchon},
  \bibinfo{author}{J.~Callen}, \bibinfo{author}{J.~Shohet},
  \bibinfo{title}{{Flux Coordinates and Magnetic Field Structure}}, Springer
  Series in Computational Physics, \bibinfo{publisher}{Springer-Verlag},
  \bibinfo{year}{1991}.
\bibitem[{Scott(2001)}]{scott_2001}
\bibinfo{author}{B.~Scott},
\newblock \bibinfo{title}{{Shifted metric procedure for flux tube treatments of
  toroidal geometry: Avoiding grid deformation}},
\newblock \bibinfo{journal}{Phys. Plasmas} \bibinfo{volume}{8}
  (\bibinfo{year}{2001}) \bibinfo{pages}{447--458}.
\bibitem[{Frankel(2004)}]{Frankel}
\bibinfo{author}{T.~Frankel}, \bibinfo{title}{{The geometry of physics: an
  introduction}}, \bibinfo{publisher}{Cambridge University Press},
  \bibinfo{edition}{second} edition, \bibinfo{year}{2004}.
\bibitem[{Wiesenberger and Held(2016)}]{FELTORv3.1}
\bibinfo{author}{M.~Wiesenberger}, \bibinfo{author}{M.~Held},
  \bibinfo{title}{Feltor v3.1}, \bibinfo{howpublished}{{Zenodo
  \url{http://doi.org/10.5281/zenodo.162407}}}, \bibinfo{year}{2016}.
\bibitem[{Hairer et~al.(1993)Hairer, N{\o}rsett, and Wanner}]{hairer93}
\bibinfo{author}{E.~Hairer}, \bibinfo{author}{S.~N{\o}rsett},
  \bibinfo{author}{G.~Wanner}, \bibinfo{title}{{Solving Ordinary Differential
  Equations I, Nonstiff Problems}}, \bibinfo{publisher}{Springer-Verlag Berlin
  Heidelberg}, \bibinfo{edition}{2nd} edition, \bibinfo{year}{1993}.
\bibitem[{Cockburn et~al.(2001)Cockburn, Kanschat, Perugia, and
  Schotzau}]{Cockburn2001}
\bibinfo{author}{B.~Cockburn}, \bibinfo{author}{G.~Kanschat},
  \bibinfo{author}{I.~Perugia}, \bibinfo{author}{D.~Schotzau},
\newblock \bibinfo{title}{{Superconvergence of the local discontinuous Galerkin
  method for elliptic problems on Cartesian grids}},
\newblock \bibinfo{journal}{SIAM J. Numer. Anal.} \bibinfo{volume}{39}
  (\bibinfo{year}{2001}) \bibinfo{pages}{264--285}.
\bibitem[{Held et~al.(2016)Held, Wiesenberger, and Stegmeir}]{Held2016}
\bibinfo{author}{M.~Held}, \bibinfo{author}{M.~Wiesenberger},
  \bibinfo{author}{A.~Stegmeir},
\newblock \bibinfo{title}{Three discontinuous galerkin schemes for the
  anisotropic heat conduction equation on non-aligned grids},
\newblock \bibinfo{journal}{Comput. Phys. Commun.} \bibinfo{volume}{199}
  (\bibinfo{year}{2016}) \bibinfo{pages}{29--39}.
\bibitem[{Einkemmer and Wiesenberger(2014)}]{Einkemmer2014}
\bibinfo{author}{L.~Einkemmer}, \bibinfo{author}{M.~Wiesenberger},
\newblock \bibinfo{title}{A conservative discontinuous galerkin scheme for the
  2d incompressible navier--stokes equations},
\newblock \bibinfo{journal}{Comput. Phys. Commun.} \bibinfo{volume}{185}
  (\bibinfo{year}{2014}) \bibinfo{pages}{2865--2873}.
\bibitem[{Cerfon and Freidberg(2010)}]{Cerfon2010}
\bibinfo{author}{A.~J. Cerfon}, \bibinfo{author}{J.~P. Freidberg},
\newblock \bibinfo{title}{"one size fits all" analytic solutions to the
  grad-shafranov equation},
\newblock \bibinfo{journal}{Phys. Plasmas} \bibinfo{volume}{17}
  (\bibinfo{year}{2010}) \bibinfo{pages}{032502}.
\bibitem[{Arnold et~al.(2002)Arnold, Brezzi, Cockburn, and Marini}]{Arnold2001}
\bibinfo{author}{D.~N. Arnold}, \bibinfo{author}{F.~Brezzi},
  \bibinfo{author}{B.~Cockburn}, \bibinfo{author}{L.~D. Marini},
\newblock \bibinfo{title}{Unified analysis of discontinuous galerkin methods
  for elliptic problems},
\newblock \bibinfo{journal}{SIAM J. Numer. Anal.} \bibinfo{volume}{39}
  (\bibinfo{year}{2002}) \bibinfo{pages}{1749--1779}.

\end{thebibliography}
\bibliographystyle{model1-num-names}


\end{document}